\documentclass[11pt]{amsart}
\usepackage[T1]{fontenc}
\usepackage{lmodern}
\usepackage{microtype}
\usepackage{mathtools,amssymb,amsthm}
\usepackage{amscd}
\usepackage[margin=1.15in]{geometry}
\usepackage[colorlinks=true,linkcolor=blue,citecolor=blue,urlcolor=blue]{hyperref}

\newtheorem{theorem}{Theorem}[section]
\newtheorem{proposition}[theorem]{Proposition}
\newtheorem{lemma}[theorem]{Lemma}
\newtheorem{corollary}[theorem]{Corollary}
\theoremstyle{definition}
\newtheorem{definition}[theorem]{Definition}
\newtheorem{example}[theorem]{Example}

\newcommand{\kHom}{\operatorname{Hom}_k}
\newcommand{\HH}{\operatorname{HH}}
\newcommand{\Ext}{\operatorname{Ext}}
\newcommand{\gr}{\operatorname{gr}}
\newcommand{\ad}{\operatorname{ad}}
\newcommand{\id}{\operatorname{id}}
\newcommand{\cD}{\mathcal D}
\newcommand{\Barres}{\mathrm B}

\title[Relative higher intertwiners]{Relative higher intertwiners and symbol calculus for quantum Borels}
\author{Abdellah Sebbar}
\address{Department of Mathematics and Statistics, University of Ottawa, Ottawa, Ontario, Canada}
\email{asebbar@uottawa.ca}
\subjclass[2020]{16T05, 16E40, 16S32, 17B37, 81R50}
\keywords{Hopf algebra, intertwiner, augmentation ideal, Hochschild cohomology, Ext group, quantum group}

\begin{document}
\begin{abstract}
Let $H$ be a Hopf algebra and let $M,N$ be left $H$-modules.  The
internal-Hom action turns $x\cdot f-\varepsilon(x)f$ into a natural measure
of the failure of a linear map $f:M\to N$ to be equivariant.  We use
iterates of this defect to define an order filtration and prove that it is
preserved by composition.  A rigidity theorem shows that the filtration
collapses to $\operatorname{Hom}_H(M,N)$ when the augmentation ideal is
idempotent.  This occurs for enveloping algebras of perfect Lie algebras
and for Drinfeld--Jimbo quantum groups when $q_i^2\ne1$.

To retain the root-direction information lost in the absolute filtration,
we introduce relative higher intertwiners for a split Hopf projection
$H\to K$.  The Rees construction connects the category of finite-order
relative higher intertwiners with its associated graded category, whose
morphisms admit a faithful realization by principal symbols.  Symbol
composition is governed by the graded coproduct.  An exact sequence identifies the
obstruction to lifting a formal symbol with a class in $\Ext_H^1$.  For a
bosonization $H=R\#K$,
relative order agrees with the braided augmentation order of $R$.

For explicit quantum-Borel modules with $q$ not a root of unity, we determine
the complete filtered endomorphism algebras.  In every simply laced finite
type, an acyclic orientation gives a quantum-space module with finite-order
endomorphism algebra
$k[X_1,\ldots,X_\ell]$; the order of $X_i$ is the number of directed paths
ending at $i$.  We also give a self-contained bar-cochain model for
$\Ext_H^\bullet(M,N)$ with internal-Hom coefficients, without
finite-dimensionality assumptions.  An explicit Hopf Alexander--Whitney composition
realizes the Yoneda product and respects the absolute order filtration,
while adjoint coefficients recover the derived center.
\end{abstract}

\maketitle

\section{Introduction}

Let $H$ be a Hopf algebra over a field $k$, with coproduct $\Delta$,
counit $\varepsilon$, and bijective antipode $S$.  If $M$ and $N$ are left
$H$-modules, then $\kHom(M,N)$ carries the familiar left action
\begin{equation}\label{eq:hom-action}
 (x\cdot f)(m)=\sum_{(x)}x_{(1)}f(S(x_{(2)})m).
\end{equation}
The invariants for this action are precisely the $H$-linear maps.  The
equivariance defect of $f$ is the map
\begin{equation}\label{eq:defect-intro}
 \langle x,f\rangle=x\cdot f-\varepsilon(x)f
\end{equation}
from $H$ to $\kHom(M,N)$.

The screening-theoretic antecedent is the Lie-de Rham construction of
Ginzburg and Schechtman \cite{GinzburgSchechtman}.  It retains the
commutator defect before integration and produces higher Ext classes from
local-system homology.  The author's subsequent constructions
\cite{Sebbar1998,SebbarCMP} develop quantum versions of this mechanism.
For quantized enveloping algebras of finite type, operator-valued
$q$-differential forms give a double complex and maps from twisted
$q$-de Rham homology to Ext groups \cite{Sebbar1998}.  For
$U_q(\widehat{\mathfrak{sl}}_2)$, quantum screening operators and
Wakimoto-type modules give the corresponding affine construction
\cite{SebbarCMP}.

The defect \eqref{eq:defect-intro} is the degree-zero bar differential
with internal-Hom coefficients.  For an enveloping algebra, its
restriction to a primitive element $X$ is the commutator
$X_Nf-fX_M$.  The present paper studies finite order under repeated
application of this defect, rather than introducing the
screening-to-Ext construction.  The composition identity already appears
in \cite{Sebbar1998}; here it is used to construct absolute and relative
filtrations.  Cohomological degree in the screening construction and
order in these filtrations are different notions, as the collapse
results below demonstrate.  We use the conventions and general framework
for Hopf actions in \cite{Montgomery}; for related structural questions
concerning Hopf filtrations, see also \cite{MontgomeryFiltrations}.

The argument follows a single structural line.  Finite equivariance order
leads to a filtered category; a split Hopf projection supplies a relative
version with nontrivial root directions; principal symbols identify its
graded structure; and bosonization makes that structure computable for
quantum Borels.  The bar construction then gives a cochain-level counterpart
of the composition law.  If
$I=\ker\varepsilon$, define
\[
 \cD^r_H(M,N)=\{f\in\kHom(M,N):I^{r+1}\cdot f=0\}.
\]
The degree-zero term is $\operatorname{Hom}_H(M,N)$.  Orders add under
composition, so the spaces
$\cD^\infty_H(M,N)$ form a filtered category.  The mechanism is the
compatibility
\[
 \Delta(I^n)\subseteq\sum_{a+b\ge n}I^a\otimes I^b,
\]
which controls every term in the Hopf action on a composite.

The absolute construction can collapse.  For
$H=U(\mathfrak g)$, the defect operators are commutators and
$I/I^2\simeq\mathfrak g/[\mathfrak g,\mathfrak g]$.  Hence the filtration
collapses for perfect Lie algebras.  A similar phenomenon occurs for
Drinfeld--Jimbo quantum groups: when $q_i^2\ne1$, the defining relations
imply $I=I^2$, so every finite-order map is already an intertwiner.  The
full augmentation ideal therefore does not retain the root-direction
information present in screening equations.  This leads to the relative
construction below.

Let $K\subset H$ be a split Hopf pair, with retraction $H\to K$ and kernel
$J$.  We require maps to be $K$-linear and measure their remaining defect
by the powers of $J$.  Relative orders again add under composition, but the
resulting filtration need not collapse.  Its associated graded category
has faithful principal symbols, with composition given by convolution
against the graded coproduct.  The Rees category provides a flat
$k[t]$-linear family whose specializations are the finite-order maps and
their symbols.
Relative higher intertwiners are also represented by explicit cyclic
$H$-modules, yielding an exact sequence in which the obstruction to
lifting a formal symbol is a connecting class in $\Ext_H^1$.  The image of
the symbol map is the kernel of this connecting map.

Radford's projection theorem identifies the corresponding filtered Hopf
algebra.  If $H=R\#K$, where $R$ is a braided Hopf algebra in a
Yetter--Drinfeld category, then
\[
 J^n=(R^+)^n\#K,
 \qquad
 \gr_JH\simeq(\gr_{R^+}R)\#K.
\]
Relative equivariance order is the braided augmentation order of
$R$.  For quantum-Borel subalgebras, this identifies the relevant graded
object with the root-height filtration of the positive or negative part.
When this part is a Nichols algebra, its braided relations, including the
mixed quantum Serre relations, remain visible in the symbol calculus.

For the quantum-Borel modules constructed in this paper, the relative
filtration has nonzero terms in every degree.  Simple-root strings produce a
polynomial endomorphism algebra
in one variable, with polynomial degree equal to relative order.  Mutually
orthogonal roots give the corresponding total-degree filtration in several
variables.  The principal calculation treats adjacent roots in every
simply laced finite type.  An acyclic orientation of the Dynkin diagram
defines an oriented quantum-space module, and every finite-order
Cartan-equivariant endomorphism of this module is polynomial.  If $X_i$
corresponds to a vertex $i$, its order is the number of directed paths
ending at $i$, including the trivial path.  Hence
\[
 \cD^\infty_{H\mid K}(V,V)\simeq k[X_1,\ldots,X_\ell]
\]
as a filtered algebra, with a weighted degree determined by the oriented
Dynkin diagram.  The weights describe the propagation of adjacent root
defects.

The same defect has a homological role.  The standard bar-resolution
identification gives
\begin{equation}\label{eq:ext-intro}
 H^n\!\left(H,\kHom(M,N)\right)\cong\Ext_H^n(M,N)
\end{equation}
for arbitrary modules $M$ and $N$ \cite{CartanEilenberg,Montgomery,Weibel}.
We include the proof to fix conventions and to avoid any use of
$M^*\otimes N$, which is unsuitable for the infinite-dimensional modules
that motivated the construction.  An explicit Hopf Alexander--Whitney
formula for cochain composition is strictly associative and induces Yoneda
composition.  The cochain filtration used for this statement is the
absolute filtration $\cD_H^r$, not the relative filtration $\cD_{H\mid K}^r$
of the quantum-Borel sections.  Thus the cohomological construction is
logically independent of the relative symbol calculus; in particular, no
relative cochain filtration is asserted here.  With adjoint coefficients the
bar complex recovers
$\HH^\bullet(H,H)\cong\Ext_H^\bullet(k,H_{\ad})$, providing the link with
derived centers of quantum groups \cite{LachowskaQi,BezrukavnikovLachowska}.

The Radford projection, the bar model for Ext, and the adjoint-coefficient
description of Hochschild cohomology are standard ingredients.  The
screening origin belongs to the earlier Lie-de Rham and $q$-de Rham
constructions \cite{GinzburgSchechtman,Sebbar1998,SebbarCMP}.  The results
specific to this paper are the relative order category and its Rees
deformation, the faithful convolution calculus for principal symbols, the
$\Ext_H^1$ obstruction that characterizes liftable symbols, and the complete
filtered endomorphism calculations for the quantum-Borel modules constructed
below.  These calculations show that the symbol theory remains nontrivial
even when the absolute augmentation filtration collapses.

The absolute group-algebra filtration belongs to the established theory of
higher order invariants.  If $A=k\Gamma$, $I$ is its augmentation ideal,
and $V$ is an $A$-module, Deitmar defines invariants of order $q$ by
\[
 H_q^0(\Gamma,V)=\{v\in V:I^qv=0\}
 \]
and studies their order-lowering maps and derived functors
\cite{DeitmarHigher,DeitmarCompact}.  Taking
$V=\kHom(M,N)$ with the conjugation action gives
\[
 \cD_H^r(M,N)=H_{r+1}^0(\Gamma,\kHom(M,N)).
\]
On the graded quotient, Deitmar's order-lowering map is induced by
\[
 [v]\longmapsto\bigl(\gamma\longmapsto(\gamma-1)v\bigr),
\]
and its iterates give injections into tensor powers of spaces of additive
characters under the hypotheses imposed there.  Thus the condition defining
absolute order, as well as the idea of recording its top difference, is not
new in the group case.  Theorem~\ref{thm:faithful-symbol} is a relative Hopf
analogue: instead of iterating first-order group differences, it evaluates
the top defect on the complete layer $J^r/J^{r+1}$ and makes this evaluation
compatible with composition through the graded coproduct.  Theorem
\ref{thm:symbol-obstruction} then determines its image by an explicit
$\Ext_H^1$ obstruction.  These features, together with the internal-Hom
interpretation, the relative Rees category, and the quantum-Borel
classifications, do not follow from the absolute higher-invariant filtration.

Passi's polynomial maps provide a second comparison.  For a group algebra
they are characterized by powers of the augmentation ideal
\cite[Chapter~V]{Passi}.  With the right-difference convention stated in
Proposition~\ref{prop:group}, absolute order is the corresponding orbit-map
condition for a linear operator under conjugation.  The present construction
differs in that the filtered objects are maps between two modules and their
composition is part of the structure.  It is also distinct from the
differential operators of Lunts and Rosenberg
\cite{LuntsRosenberg}.  Their order is defined by commutators with
multiplication in a chosen algebra, whereas our order is defined by a
symmetry Hopf algebra and requires no module-algebra structure on either
$M$ or $N$.  The cyclic modules $A_r$ in
\eqref{eq:truncated-relative-module} play the role of relative principal
parts: they represent operators of order at most $r$, and their successive
layers carry the principal symbols.  This representation, together with
the obstruction sequence and its bosonization description, is the main
structural difference from the preceding operator filtrations.

Section 2 establishes the internal-Hom
and defect identities, and Section 3 proves filtered composition.  Section
4 identifies the collapse of absolute order in the principal classical
and quantum examples.  Section 5 develops relative order, the Rees and
symbol constructions, the lifting obstruction, bosonization, and the
quantum-Borel classification.  Section 6 gives the standard Ext comparison
and constructs the Hopf Alexander--Whitney composition together with its
absolute filtration.  Its final
subsection treats adjoint coefficients and the derived center.  Section 7
compares Lie-de Rham and quantum screening cocycles with the order
filtrations considered here.

\section{Equivariance defects and the categorical problem}

An order theory for equivariance defects first requires a composition law.
This section establishes that law in the internal-Hom formalism.
Throughout, tensor products and linear maps are taken over $k$.  We assume
that $S$ is bijective.  The latter hypothesis is automatic for the principal
examples considered below, and it allows symmetric changes between the
regular and diagonal actions.

For left $H$-modules $M,N$, define the action \eqref{eq:hom-action}.  A
direct calculation using the Hopf identities shows that this is a left
$H$-module structure.

\begin{proposition}\label{prop:adjunction}
For left $H$-modules $L,M,N$, there is a natural isomorphism
\[
 \operatorname{Hom}_H(L\otimes M,N)
 \cong \operatorname{Hom}_H(L,\kHom(M,N)),
\]
where $L\otimes M$ has the diagonal action and the internal Hom has the
action \eqref{eq:hom-action}.
\end{proposition}

\begin{proof}
For a linear map $F:L\otimes M\to N$, let $\widetilde F(\ell)(m)=F(\ell
\otimes m)$.  Suppose first that $F$ is $H$-linear.  For $x\in H$,
\begin{align*}
 (x\cdot\widetilde F(\ell))(m)
 &=\sum x_{(1)}F(\ell\otimes S(x_{(2)})m)\\
 &=\sum F(x_{(1)}\ell\otimes x_{(2)}S(x_{(3)})m)\\
 &=F(x\ell\otimes m)=\widetilde F(x\ell)(m).
\end{align*}
Thus $\widetilde F$ is $H$-linear.  Conversely, if $\widetilde F$ is
$H$-linear, then
\begin{align*}
 F(x_{(1)}\ell\otimes x_{(2)}m)
 &=(x_{(1)}\cdot\widetilde F(\ell))(x_{(2)}m)\\
 &=\sum x_{(1)}\widetilde F(\ell)(S(x_{(2)})x_{(3)}m)\\
 &=xF(\ell\otimes m).
\end{align*}
The two assignments are inverse and natural.
\end{proof}

\begin{proposition}\label{prop:composition}
For left $H$-modules $M,N,P$, the composition map
\[
 \kHom(N,P)\otimes\kHom(M,N)\longrightarrow\kHom(M,P)
\]
is $H$-linear.  Equivalently,
\begin{equation}\label{eq:composition}
 x\cdot(f\circ g)=\sum_{(x)}(x_{(1)}\cdot f)\circ(x_{(2)}\cdot g).
\end{equation}
\end{proposition}

\begin{proof}
For $m\in M$, the right side of \eqref{eq:composition}, evaluated at $m$,
is
\[
 \sum x_{(1)}f\bigl(S(x_{(2)})x_{(3)}g(S(x_{(4)})m)\bigr)
 =\sum x_{(1)}f\bigl(g(S(x_{(2)})m)\bigr),
\]
by coassociativity and the antipode identity.  This is the left side.
No dualizability or finite-dimensionality assumption occurs.
\end{proof}

For $x\in H$ and $f\in\kHom(M,N)$, put
\[
 \langle x,f\rangle=x\cdot f-\varepsilon(x)f.
\]
The following formulas will be used repeatedly.

\begin{proposition}\label{prop:defect}
For $x,y\in H$ and composable maps $f,g$,
\begin{align}
 \langle xy,f\rangle
 &=x\cdot\langle y,f\rangle+\varepsilon(y)\langle x,f\rangle,
 \label{eq:defect-product}\\
 \langle x,f\circ g\rangle
 &=\sum_{(x)}\langle x_{(1)},f\rangle\circ
       \langle x_{(2)},g\rangle
   +\langle x,f\rangle\circ g+f\circ\langle x,g\rangle.
 \label{eq:defect-composition}
\end{align}
In particular, $f$ is $H$-linear if and only if
$\langle x,f\rangle=0$ for every $x\in H$.
\end{proposition}

\begin{proof}
Since the action of $H$ on the internal Hom is a representation,
\begin{align*}
 \langle xy,f\rangle
 &=x\cdot(y\cdot f)-\varepsilon(x)\varepsilon(y)f\\
 &=x\cdot\bigl(y\cdot f-\varepsilon(y)f\bigr)
   +\varepsilon(y)\bigl(x\cdot f-\varepsilon(x)f\bigr),
\end{align*}
which is \eqref{eq:defect-product}.  For \eqref{eq:defect-composition},
insert
\[
 x_{(i)}\cdot f
 =\langle x_{(i)},f\rangle+\varepsilon(x_{(i)})f
\]
and the analogous formula for $g$ into \eqref{eq:composition}.  Expanding
the product gives four terms.  The term containing two counits is
$\varepsilon(x)f\circ g$.  Applying
$\sum\varepsilon(x_{(1)})x_{(2)}=x$ and
$\sum x_{(1)}\varepsilon(x_{(2)})=x$ to the two mixed terms gives
$f\circ\langle x,g\rangle$ and
$\langle x,f\rangle\circ g$, respectively.  Subtracting
$\varepsilon(x)f\circ g$ gives \eqref{eq:defect-composition}.

If $f$ is $H$-linear, then \eqref{eq:hom-action} gives
$x\cdot f=\varepsilon(x)f$.  Conversely, this equality implies
$xf(m)=f(xm)$ after applying the internal-Hom adjunction of
Proposition~\ref{prop:adjunction}.  Hence the defect vanishes exactly on
$H$-linear maps.
\end{proof}

\section{Absolute order and composability}

The defect identities define an order calculus.  Let
$I=\ker\varepsilon$ be the augmentation ideal.  We use the convention
$I^0=H$.

\begin{definition}\label{def:higher-intertwiner}
For $r\ge0$, the space of \emph{higher $H$-intertwiners of order at most
$r$} from $M$ to $N$ is
\[
 \cD^r_H(M,N)=\{f\in\kHom(M,N):I^{r+1}\cdot f=0\}.
\]
We also put $\cD^{-1}_H(M,N)=0$ and
$\cD^\infty_H(M,N)=\bigcup_{r\ge0}\cD^r_H(M,N)$.
\end{definition}

For a nonzero $f\in\cD^\infty_H(M,N)$, its order is the least $r$ for
which $f\in\cD^r_H(M,N)$.  We set $\operatorname{ord}(0)=-\infty$.

Since $I\cdot f=0$ is equivalent to $x\cdot f=\varepsilon(x)f$ for all
$x\in H$, one has
\begin{equation}\label{eq:order-zero}
 \cD^0_H(M,N)=\operatorname{Hom}_H(M,N).
\end{equation}

\begin{lemma}\label{lem:filtered-coproduct}
For every $n\ge0$,
\begin{equation}\label{eq:filtered-coproduct}
 \Delta(I^n)\subseteq
 \sum_{a+b\ge n} I^a\otimes I^b.
\end{equation}
Here the sum may be restricted to $0\le a,b\le n$ after replacing terms
of higher exponent by smaller powers.
\end{lemma}

\begin{proof}
The counit identities imply
$\Delta(I)\subseteq I\otimes H+H\otimes I$.  Indeed, the image of
$\Delta(x)$ in $(H/I)\otimes(H/I)\cong k$ is $\varepsilon(x)$, which is
zero for $x\in I$.  Since $\Delta$ is an algebra homomorphism,
\[
 \Delta(I^n)\subseteq(I\otimes H+H\otimes I)^n.
\]
Every product in the expansion belongs to $I^a\otimes I^b$ with
$a+b=n$.  Allowing larger exponents gives \eqref{eq:filtered-coproduct}.
\end{proof}

\begin{theorem}\label{thm:main}
Let $M,N,P$ be left $H$-modules.  For all $p,q\ge0$,
\begin{equation}\label{eq:main}
 \cD^p_H(N,P)\circ\cD^q_H(M,N)
 \subseteq\cD^{p+q}_H(M,P).
\end{equation}
Consequently, the objects consisting of left $H$-modules and the morphism
spaces $\cD^\infty_H(M,N)$ form a filtered $k$-linear category.  Its
degree-zero subcategory is $H\text{-}\mathrm{Mod}$.
\end{theorem}

\begin{proof}
Take $f\in\cD^p_H(N,P)$ and $g\in\cD^q_H(M,N)$.  If
$x\in I^{p+q+1}$, Lemma~\ref{lem:filtered-coproduct} writes $\Delta(x)$ as
a sum of terms $u\otimes v$ with $u\in I^a$, $v\in I^b$, and
$a+b\ge p+q+1$.  For every such pair, either $a\ge p+1$ or $b\ge q+1$.
By Proposition~\ref{prop:composition}, each corresponding summand of
$x\cdot(f\circ g)$ is $(u\cdot f)\circ(v\cdot g)$ and therefore vanishes.
Thus $I^{p+q+1}\cdot(f\circ g)=0$.

The identity map of an $H$-module is $H$-linear and hence has order zero.
Associativity is inherited from ordinary composition.  The filtration law
is \eqref{eq:main}, and \eqref{eq:order-zero} identifies its degree-zero
subcategory.
\end{proof}

\begin{corollary}\label{cor:graded}
There is a graded $k$-linear category with the same objects as
$H\text{-}\mathrm{Mod}$ and graded morphisms
\[
 \gr^r\cD_H(M,N)=\cD^r_H(M,N)/\cD^{r-1}_H(M,N).
\]
Composition is induced by ordinary composition and has degree addition.
Moreover,
\[
 I^{r+1}\cdot\gr^r\cD_H(M,N)=0.
\]
\end{corollary}

\begin{proof}
The spaces in the quotient are $H$-submodules.  Indeed, if
$f\in\cD_H^r(M,N)$, $x\in H$, and $a\in I^{r+1}$, then
\[
 a\cdot(x\cdot f)=(ax)\cdot f=0,
\]
because $I^{r+1}$ is a right ideal.  Thus the action on
$\gr^r\cD_H(M,N)$ is well defined.

If $f\in\cD^p_H(N,P)$ and $g\in\cD^q_H(M,N)$, then
Theorem~\ref{thm:main} places $f\circ g$ in
$\cD^{p+q}_H(M,P)$.  If $f$ is changed by an element of
$\cD^{p-1}_H(N,P)$, the composite changes by an element of
$\cD^{p+q-1}_H(M,P)$; the same holds when $g$ is changed by an element of
$\cD^{q-1}_H(M,N)$.  Composition therefore descends to the indicated
quotients.  For a representative $f\in\cD^r_H(M,N)$, the defining
condition is $I^{r+1}\cdot f=0$, so the induced action of $I^{r+1}$ on its
class in $\gr^r\cD_H(M,N)$ is zero.
\end{proof}

The filtration is functorial.

\begin{proposition}\label{prop:functoriality}
Let $a:M'\to M$ and $b:N\to N'$ be $H$-linear.  Then
\[
 f\longmapsto b\circ f\circ a
\]
maps $\cD^r_H(M,N)$ into $\cD^r_H(M',N')$.  In particular, the
construction descends to full subcategories and is stable under
isomorphisms of modules.
\end{proposition}

\begin{proof}
Both $a$ and $b$ are $H$-linear, so they have order zero by
\eqref{eq:order-zero}.  If $f$ has order at most $r$, then
Theorem~\ref{thm:main} gives
\[
 f\circ a\in\cD^r_H(M',N),
 \qquad
 b\circ f\circ a\in\cD^r_H(M',N').
\]
\end{proof}

\begin{proposition}\label{prop:tensor}
Assume that $H$ is cocommutative.  If $f\in\cD^p_H(M,N)$ and
$g\in\cD^q_H(M',N')$, then the ordinary tensor map
$f\otimes g:M\otimes M'\to N\otimes N'$ has order at most $p+q$.
\end{proposition}

\begin{proof}
For $m\in M$ and $m'\in M'$, the internal-Hom action gives
\begin{align*}
 \bigl(x\cdot(f\otimes g)\bigr)(m\otimes m')
 &=\sum_{(x)}x_{(1)}f(S(x_{(4)})m)\otimes
                  x_{(2)}g(S(x_{(3)})m')\\
 &=\sum_{(x)}x_{(1)}f(S(x_{(2)})m)\otimes
                  x_{(3)}g(S(x_{(4)})m').
\end{align*}
The second equality uses cocommutativity to permute the middle tensor
factors in the iterated coproduct.  The last expression is
\[
 \sum_{(x)}\bigl((x_{(1)}\cdot f)\otimes
 (x_{(2)}\cdot g)\bigr)(m\otimes m').
\]
Let $x\in I^{p+q+1}$.  By Lemma~\ref{lem:filtered-coproduct}, each tensor
in an expansion of $\Delta(x)$ has one factor in $I^{p+1}$ or the other in
$I^{q+1}$.  The corresponding summand above is zero because
$I^{p+1}\cdot f=0$ and $I^{q+1}\cdot g=0$.  Hence
$I^{p+q+1}\cdot(f\otimes g)=0$.
\end{proof}

For a general Hopf algebra, the covariance identity used in this proof
need not hold for the ordinary tensor product of linear maps.  It requires
an interchange of the two middle factors in the corresponding internal-Hom
objects, and that interchange is not $H$-linear in an unbraided module
category.  We therefore make no tensor-product closure assertion outside
the cocommutative case.  Under suitable quasitriangular hypotheses one may
instead use the braided tensor product, which is a different construction.

The filtration need not be exhaustive.  Its union consists of the maps on
which the augmentation ideal acts locally nilpotently.  In particular, it
is exhaustive if the $I$-action on $\kHom(M,N)$ is locally nilpotent.

\section{The collapse problem for absolute order}

The absolute filtration carries information only if it has positive layers.
This section determines when those layers disappear.
The enveloping and group-algebra calculations identify the classical
commutator and finite-difference models; the quantum-group calculation
then gives the obstruction that forces the relative construction of
Section~5.

\subsection{Universal enveloping algebras}

Let $H=U(\mathfrak g)$ with its standard cocommutative Hopf structure.
For $X\in\mathfrak g$ and $f\in\kHom(M,N)$,
\begin{equation}\label{eq:lie-defect}
 \langle X,f\rangle(m)=Xf(m)-f(Xm).
\end{equation}
Since $I$ is generated by $\mathfrak g$, Definition
\ref{def:higher-intertwiner} has the following concrete form.

\begin{proposition}\label{prop:lie}
A map $f:M\to N$ belongs to $\cD^r_{U(\mathfrak g)}(M,N)$ if and only if
\[
 [X_0,[X_1,\ldots,[X_r,f]\ldots]]=0
\]
for all $X_0,\ldots,X_r\in\mathfrak g$, where the commutator is interpreted
as in \eqref{eq:lie-defect}.  The composition of maps of commutator orders
$p$ and $q$ has commutator order at most $p+q$.
\end{proposition}

\begin{proof}
The iterated commutator is exactly
$(X_0\cdots X_r)\cdot f$.  Since $I$ is generated as a nonunital algebra
by $\mathfrak g$,
$I^{r+1}$ is spanned by products of at least $r+1$ elements of
$\mathfrak g$.  The first assertion follows, and the second is
Theorem~\ref{thm:main}.
\end{proof}

\begin{proposition}
\label{prop:indecomp-Ug}
For every Lie algebra $\mathfrak g$ there is a natural isomorphism
\[
 I/I^2\cong \mathfrak g/[\mathfrak g,\mathfrak g].
\]
Consequently, if $\mathfrak g$ is perfect, then $I=I^2$ and
\[
 \cD^r_{U(\mathfrak g)}(M,N)
 =\operatorname{Hom}_{U(\mathfrak g)}(M,N)
 \qquad(r\ge0).
\]
\end{proposition}

\begin{proof}
The inclusion $\mathfrak g\to I$ induces a surjection onto $I/I^2$,
since every product of two or more elements of $\mathfrak g$ lies in
$I^2$.  The relation $[X,Y]=XY-YX$ shows that
$[\mathfrak g,\mathfrak g]$ maps to zero.  We therefore obtain a
surjective map
\[
 \alpha:\mathfrak g/[\mathfrak g,\mathfrak g]\longrightarrow I/I^2.
\]
For injectivity, let $\mathfrak g_{\mathrm{ab}}$ denote the abelianization
and form the algebra $k\oplus\mathfrak g_{\mathrm{ab}}$ in which the second
summand has square zero.  The map $X\mapsto(0,\bar X)$ is a Lie algebra
homomorphism from $\mathfrak g$ to this algebra with the commutator bracket.
By the universal property of $U(\mathfrak g)$, it extends to an algebra
homomorphism
\[
 U(\mathfrak g)\longrightarrow k\oplus\mathfrak g_{\mathrm{ab}}.
\]
Its restriction to $I$ vanishes on $I^2$ and induces a map
$I/I^2\to\mathfrak g_{\mathrm{ab}}$ inverse to $\alpha$.
If $\mathfrak g$ is perfect, then $I/I^2=0$, hence $I=I^2$, and the last
assertion follows from Definition~\ref{def:higher-intertwiner}.
\end{proof}

\begin{example}
Let $\mathfrak g=kX$ act on $k[t]$ by $X=d/dt$.  Multiplication by a
polynomial $a(t)$ satisfies
\[
 [X,m_a]=m_{a'},\qquad
 [X,[X,\ldots,[X,m_a]\ldots]]=m_{a^{(s)}}.
\]
Hence $m_a$ has order at most $\deg a$.  Theorem~\ref{thm:main} recovers
$\deg(ab)\le\deg a+\deg b$ as the composition law for multiplication
operators.  More generally, the filtration is a representation-theoretic
version of the commutator definition of differential operators.
\end{example}

\subsection{Group algebras and finite differences}

Let $H=kG$.  For $s\in G$,
\[
 (s\cdot f)(m)=s f(s^{-1}m),\qquad
 \langle s,f\rangle=sfs^{-1}-f.
\]
Since $I$ is spanned by $s-1$, the order condition becomes an iterated
finite-difference condition.

\begin{proposition}\label{prop:group}
For $f\in\kHom(M,N)$, the following are equivalent:
\begin{enumerate}
\item $f\in\cD^r_{kG}(M,N)$;
\item for every $s_0,\ldots,s_r\in G$,
\[
 (s_0-1)\cdots(s_r-1)\cdot f=0;
\]
\item every $(r+1)$-fold finite difference of the orbit map
$s\mapsto sfs^{-1}$ vanishes at the identity.
\end{enumerate}
\end{proposition}

\begin{proof}
The augmentation ideal of $kG$ is generated as a vector space and as an
ideal by the elements $s-1$.  Thus $I^{r+1}$ is spanned by their products.
It follows that conditions (1) and (2) are equivalent.  For the orbit map
$F(s)=sfs^{-1}$, define
\[
 (\delta_uF)(s)=F(su)-F(s).
\]
Evaluation at the identity gives
$(\delta_uF)(1)=(u-1)\cdot f$.  Iterating this identity identifies an
$(r+1)$-fold difference at the identity with the action of a product
$(s_0-1)\cdots(s_r-1)$ on $f$.  Hence (2) and (3) are equivalent.
\end{proof}

For an arbitrary $kG$-module $V$, the spaces annihilated by $I^{r+1}$ are
the higher order invariants studied by Deitmar
\cite{DeitmarHigher,DeitmarCompact}.  Proposition~\ref{prop:group} is their
specialization to the coefficient module $V=\kHom(M,N)$ with the
conjugation action.  The order-lowering map in that theory sends the class
of $v$ to the function $s\mapsto(s-1)v$; here this is the first equivariance
defect $s\mapsto sfs^{-1}-f$.

With the right-difference convention used above, the same condition is the
orbit-map version of Passi's augmentation-ideal criterion for polynomial
maps \cite[Chapter~V]{Passi}.  The additional structure here is that $F$ is
the orbit map of an operator, so composition of operators induces the
filtered composition studied in Section~3.  In the relative construction
of Section~5, one first imposes equivariance for the retract $K$ and then
uses powers of the projection kernel $J$.  This is different from imposing
the full augmentation powers on a single group module.

If $G$ is finite and $\operatorname{char}k\nmid |G|$, the augmentation
ideal is idempotent.  Hence the filtration collapses:
\[
 \cD^r_{kG}(M,N)=\operatorname{Hom}_{kG}(M,N)\qquad(r\ge0).
\]
Indeed, the central idempotent
$1-|G|^{-1}\sum_{g\in G}g$ generates $I$.  Positive-order terms can
therefore occur only outside this semisimple finite-group case.

\subsection{Quantized enveloping algebras}

In this subsection, $\mathfrak g$ is a finite-dimensional semisimple Lie
algebra with finite-type symmetrizable Cartan matrix $(a_{ij})$.  We fix
an element $q\in k^\times$ and positive integers $d_i$ such that
$d_i a_{ij}=d_j a_{ji}$, put
$q_i=q^{d_i}$, and use the simply connected Drinfeld--Jimbo presentation.
Whenever the relation below is used, we assume $q_i-q_i^{-1}\ne0$.

Let $U_q(\mathfrak g)$ be the resulting quantum group, with
$\varepsilon(E_i)=\varepsilon(F_i)=0$ and $\varepsilon(K_i)=1$.  Its
augmentation ideal is generated by
\[
 E_i,\qquad F_i,\qquad K_i-1.
\]
Consequently a linear map has order at most $r$ precisely when every word
of length $r+1$ in the corresponding defect operators annihilates it.
For example,
\begin{align*}
 \langle E_i,f\rangle&=E_i\cdot f,\\
 \langle F_i,f\rangle&=F_i\cdot f,\\
 \langle K_i,f\rangle&=K_i\cdot f-f.
\end{align*}
The last expression is a multiplicative $q$-difference, while the first
two are skew derivation defects.  Relations among these defects are
inherited automatically from the defining relations of
$U_q(\mathfrak g)$.

\begin{theorem}\label{thm:quantum-collapse}
Assume $\operatorname{char}k=0$ and $q_i^2\ne1$ for every simple root.
For the Drinfeld--Jimbo algebra $U_q(\mathfrak g)$ one has $I=I^2$.
Consequently,
\[
 \cD^r_{U_q(\mathfrak g)}(M,N)
 =\operatorname{Hom}_{U_q(\mathfrak g)}(M,N)
 \qquad(r\ge0)
\]
for all modules $M,N$.
\end{theorem}

\begin{proof}
It is enough to show that $E_i,F_i,K_i-1$ lie in $I^2$.  The relation
$K_iE_i=q_i^2E_iK_i$ gives
\[
 (K_i-1)E_i-q_i^2E_i(K_i-1)=(q_i^2-1)E_i,
\]
so $E_i\in I^2$.  Similarly, $K_iF_i=q_i^{-2}F_iK_i$ gives
$F_i\in I^2$.  The relation
\[
 [E_i,F_i]=\frac{K_i-K_i^{-1}}{q_i-q_i^{-1}}
\]
then implies $K_i-K_i^{-1}\in I^2$.  Modulo $I^2$ one has
$K_i^{-1}\equiv1-(K_i-1)$, and hence
$K_i-K_i^{-1}\equiv2(K_i-1)$.  Since $2$ is invertible,
$K_i-1\in I^2$.  Therefore $I=I^2$, and the assertion about higher
intertwiners follows.
\end{proof}

The theorem rules out the full augmentation filtration as an order
filtration for screening operators.  Their $q$-difference equations define
cochains and Ext classes, but no positive-order higher intertwiners for
this filtration.  The relative filtration below retains the Borel and
root-direction data.

\section{Relative order, symbols, and quantum Borels}

The collapse in Theorem~\ref{thm:quantum-collapse} results from mixing
the Cartan and root directions inside the full augmentation ideal.  The
following relative construction separates them.

\begin{definition}
A \emph{split Hopf pair} consists of Hopf algebra morphisms
\[
 K\xrightarrow{\iota}H\xrightarrow{\pi}K,
 \qquad \pi\iota=\id_K.
\]
We identify $K$ with its image in $H$ and put $J=\ker\pi$.  Thus $J$ is
a Hopf ideal and $H=K\oplus J$ as vector spaces.
\end{definition}

Every $H$-module is a $K$-module by restriction.  The internal-Hom action
of $H$ restricts to the usual internal-Hom action of $K$.

\begin{definition}\label{def:relative}
For a split Hopf pair $(H,K)$ and left $H$-modules $M,N$, define
\[
 \cD^r_{H\mid K}(M,N)
 =
 \{f\in\operatorname{Hom}_K(M,N):J^{r+1}\cdot f=0\}.
\]
Its elements are the \emph{relative higher intertwiners of order at most
$r$}.
\end{definition}

For a nonzero relative higher intertwiner, write
\[
 \operatorname{ord}_{H\mid K}(f)
 =\min\{r:f\in\cD^r_{H\mid K}(M,N)\},
\qquad \operatorname{ord}_{H\mid K}(0)=-\infty.
\]
When the split pair is fixed throughout an example, we abbreviate this
order to $\operatorname{ord}(f)$.

The relative order-zero term has the expected meaning.

\begin{lemma}\label{lem:relative-zero}
For every split Hopf pair,
\[
 \cD^0_{H\mid K}(M,N)=\operatorname{Hom}_H(M,N).
\]
\end{lemma}

\begin{proof}
If $f$ is $H$-linear, it is $K$-linear and $J\cdot f=0$.  Conversely,
write $x=\iota\pi(x)+(x-\iota\pi(x))$ with the second term in $J$.
The defect of $f$ vanishes on the first summand by $K$-linearity and on
the second by $J\cdot f=0$.  Hence it vanishes on all of $H$.
\end{proof}

\begin{theorem}\label{thm:relative-composition}
Let $(H,K)$ be a split Hopf pair and let $M,N,P$ be left $H$-modules.
Then
\[
 \cD^p_{H\mid K}(N,P)\circ\cD^q_{H\mid K}(M,N)
 \subseteq\cD^{p+q}_{H\mid K}(M,P).
\]
Thus the spaces $\cD^\infty_{H\mid K}(M,N)$ form a filtered category
whose degree-zero subcategory is $H\text{-}\mathrm{Mod}$.
\end{theorem}

\begin{proof}
The composite of two $K$-linear maps is $K$-linear.  Since $J$ is a Hopf
ideal,
\[
 \Delta(J^n)\subseteq\sum_{a+b\ge n}J^a\otimes J^b.
\]
Indeed, $\Delta(J)\subseteq J\otimes H+H\otimes J$, and multiplication of
$n$ such sums gives the displayed inclusion.  Let
$f\in\cD^p_{H\mid K}(N,P)$, $g\in\cD^q_{H\mid K}(M,N)$, and
$x\in J^{p+q+1}$.  In every term $u\otimes v$ of $\Delta(x)$, either
$u\in J^{p+1}$ or $v\in J^{q+1}$.  Hence $u\cdot f=0$ or
$v\cdot g=0$.  Proposition~\ref{prop:composition} expresses
$x\cdot(f\circ g)$ as a sum of composites
$(u\cdot f)\circ(v\cdot g)$, all of which vanish.  Thus
$f\circ g\in\cD^{p+q}_{H\mid K}(M,P)$.  Lemma~\ref{lem:relative-zero}
identifies the degree-zero term.
\end{proof}

\begin{corollary}\label{cor:relative-graded}
The quotients
\[
 \gr^r\cD_{H\mid K}(M,N)
 =\cD^r_{H\mid K}(M,N)/\cD^{r-1}_{H\mid K}(M,N)
\]
are the morphism spaces of a graded category, with composition of degree
addition.
\end{corollary}

\begin{proof}
The graded composition follows from
Theorem~\ref{thm:relative-composition}, by the same quotient argument as
in Corollary~\ref{cor:graded}.
\end{proof}

The $H$-action on $\kHom(M,N)$ does not in general preserve the subspace
$\operatorname{Hom}_K(M,N)$.  Consequently the relative associated graded
morphism spaces do not carry a natural action of $\gr_JH$.  The correct
replacement records all values $x\cdot f$ simultaneously.  Their dependence
on $x$ is $K$-equivariant for the adjoint action on $J^r/J^{r+1}$, and this
gives the principal-symbol map constructed below.

The filtered category has a canonical deformation to its associated graded
category.
For either the absolute or relative filtration, write simply
$\cD^r(M,N)$ and define
\[
 \mathcal R\cD(M,N)=
 \bigoplus_{r\geq0}\cD^r(M,N)t^r
 \subseteq \kHom(M,N)[t].
\]

\begin{proposition}\label{prop:rees-category}
The spaces $\mathcal R\cD(M,N)$ are the morphism spaces of a
$k[t]$-linear category under
\[
 (f t^p)\circ(g t^q)=(f\circ g)t^{p+q}.
\]
Every morphism space is flat over $k[t]$, and there are natural
specializations
\[
 \mathcal R\cD(M,N)/(t-1)\mathcal R\cD(M,N)
 \simeq\cD^\infty(M,N),
\]
\[
 \mathcal R\cD(M,N)/t\mathcal R\cD(M,N)
 \simeq\gr\cD(M,N).
\]
Thus the Rees construction is a flat family whose fibers at $t=1$ and
$t=0$ are the category of finite-order higher intertwiners and its associated
graded category, respectively.
\end{proposition}

\begin{proof}
The composition law follows from Theorem~\ref{thm:main} in the absolute
case and Theorem~\ref{thm:relative-composition} in the relative case.
The Rees module is a submodule of the free $k[t]$-module
$\kHom(M,N)[t]$, hence is torsion-free and therefore flat over the
principal ideal domain $k[t]$.  Evaluation at $t=1$ has image
$\cD^\infty(M,N)$.  Suppose
$F(t)=\sum_{i=0}^n f_it^i$ belongs to its kernel.  Then
$\sum_{i=0}^n f_i=0$.  For $0\le s<n$, put
\[
 q_s=-\sum_{i=0}^s f_i.
\]
The filtration is increasing, so $q_s\in\cD^s(M,N)$.  Comparing
coefficients gives
\[
 F(t)=(t-1)\sum_{s=0}^{n-1}q_st^s.
\]
Hence the kernel of evaluation is
$(t-1)\mathcal R\cD(M,N)$, which proves the first specialization.  In
degree $r$, multiplication by
$t$ identifies the degree-$r$ part of $t\mathcal R\cD(M,N)$ with
$\cD^{r-1}(M,N)t^r$.  Hence the second quotient is
$\bigoplus_r\cD^r/\cD^{r-1}$, with its graded composition.
\end{proof}

\subsection{Faithful symbols and bosonization}

The associated graded category has a faithful principal-symbol
realization.  Put
$Q_r=J^r/J^{r+1}$.  The adjoint action of $K$ preserves $J^r$ and
$J^{r+1}$, so $Q_r$ is a $K$-module.  We give $Q_r\otimes M$ the diagonal
$K$-action.  For $f\in\cD^r_{H\mid K}(M,N)$ define
\[
 \sigma_r(f):Q_r\otimes M\longrightarrow N,
 \qquad \sigma_r(f)(\bar x\otimes m)=(x\cdot f)(m).
\]

We use the descending filtration $J^0=H\supseteq J\supseteq J^2\supseteq
\cdots$.  On $H\otimes H$ its total filtration is
\begin{equation}\label{eq:total-tensor-filtration}
 F^n(H\otimes H)=\sum_{a+b=n}J^a\otimes J^b
 =\sum_{a+b\geq n}J^a\otimes J^b.
\end{equation}
The two sums agree because the powers of $J$ form a descending sequence.
Since $J$ is a Hopf ideal, $\Delta(J^n)\subseteq F^n(H\otimes H)$.
Consequently the coproduct of $H$ induces
\[
 Q_n\longrightarrow\bigoplus_{a+b=n}Q_a\otimes Q_b.
\]
We denote its $(p,q)$ component by $\Delta_{p,q}$.

\begin{theorem}\label{thm:faithful-symbol}
The map above induces a natural injection
\[
 \sigma_r:\gr^r\cD_{H\mid K}(M,N)
 \lhook\joinrel\longrightarrow \operatorname{Hom}_K(Q_r\otimes M,N).
\]
Let $f$ and $g$ have orders at most $p$ and $q$, respectively.  For the
component $\Delta_{p,q}:Q_{p+q}\to Q_p\otimes Q_q$ of the graded
coproduct, one has
\begin{equation}\label{eq:symbol-composition}
 \sigma_{p+q}([f\circ g])(\bar x\otimes m)
 =\sum \sigma_p([f])\bigl(\bar x_{(1)}\otimes
       \sigma_q([g])(\bar x_{(2)}\otimes m)\bigr),
\end{equation}
where $\Delta_{p,q}(\bar x)=\sum\bar x_{(1)}\otimes\bar x_{(2)}$.
Consequently, the family $\{\sigma_r\}_{r\geq0}$ is injective on every
graded morphism space and intertwines graded composition with the convolution
operation defined by the right side of \eqref{eq:symbol-composition}.  No
surjectivity onto the full target is asserted: not every $K$-equivariant map
$Q_r\otimes M\to N$ need lift to a relative higher intertwiner.
\end{theorem}

\begin{proof}
First we check equivariance.  The curried form of the symbol sends
$\bar x$ to $x\cdot f\in\kHom(M,N)$.  For $k\in K$, the adjoint action
and the $K$-invariance of $f$ give
\[
 (k\mathbin{\triangleright}x)\cdot f
 =\sum k_{(1)}\cdot\bigl(x\cdot(S(k_{(2)})\cdot f)\bigr)
 =k\cdot(x\cdot f).
\]
Thus the curried map $Q_r\to\kHom(M,N)$ is $K$-linear, and
Proposition~\ref{prop:adjunction} shows that $\sigma_r(f)$ is $K$-linear.
If $x\in J^r$ and $y\in J$, then
$y\cdot(x\cdot f)=(yx)\cdot f=0$.  Hence the curried symbol takes values
in
\[
 V^J=\{v\in\kHom(M,N):Jv=0\}.
\]
Equivalently, the image of $\sigma_r$ lies in the subspace corresponding
to $\operatorname{Hom}_K(Q_r,V^J)$ under the internal-Hom adjunction.

If $x$ is changed by an element of $J^{r+1}$, the value is unchanged
because $J^{r+1}\cdot f=0$.  If $f$ is changed by an element of
$\cD^{r-1}$, then $J^r$ annihilates the change.  Thus $\sigma_r$ is well
defined on the displayed quotient.  If its value is zero, then
$J^r\cdot f=0$, which says precisely that $f\in\cD^{r-1}$.  This proves
injectivity.

For composition, Proposition~\ref{prop:composition} gives
\[
 x\cdot(f\circ g)=\sum_{(x)}(x_{(1)}\cdot f)\circ(x_{(2)}\cdot g).
\]
Pass to the quotient
$F^{p+q}(H\otimes H)/F^{p+q+1}(H\otimes H)$ defined in
\eqref{eq:total-tensor-filtration}.  Terms of total filtration degree
greater than $p+q$ vanish there, and the degree-$p+q$ quotient is the direct
sum of the bidegrees $(a,b)$ with $a+b=p+q$.  A component with $a>p$
annihilates $f$, while a component with $b>q$ annihilates $g$.  Since
$a+b=p+q$, every bidegree other than $(p,q)$ has one of these two
properties.  Only the $(p,q)$ component remains.  Evaluation on $m$ gives
\eqref{eq:symbol-composition}.
\end{proof}

The target in Theorem~\ref{thm:faithful-symbol} contains formal symbols
which need not lift.  The next result identifies the image and gives an
exact obstruction to lifting a formal symbol.

Put $K^+=\ker(\varepsilon|_K)$.  For $r\geq0$, define the cyclic left
$H$-module
\begin{equation}\label{eq:truncated-relative-module}
 A_r=H/(HK^++J^{r+1}),
\end{equation}
and, for $r\geq1$, let
\begin{equation}\label{eq:relative-layer}
 C_r=(HK^++J^r)/(HK^++J^{r+1}).
\end{equation}
There is a natural short exact sequence
\begin{equation}\label{eq:relative-layer-sequence}
 0\longrightarrow C_r\longrightarrow A_r
 \longrightarrow A_{r-1}\longrightarrow0.
\end{equation}
The quotient map $J^r\to C_r$ factors through a natural surjection
\begin{equation}\label{eq:Qr-Cr}
 \rho_r:Q_r\longrightarrow C_r,
 \qquad x+J^{r+1}\longmapsto x+(HK^++J^{r+1}).
\end{equation}
Its kernel is
\begin{equation}\label{eq:rho-kernel}
 \ker\rho_r=
 \frac{J^r\cap(HK^++J^{r+1})}{J^{r+1}}.
\end{equation}
Here $Q_r$ has the adjoint $K$-action and $C_r$ has the action induced by
left multiplication.  The map $\rho_r$ is $K$-linear.  Indeed, for
$k\in K$ and $x\in J^r$,
\[
 \rho_r(k\mathbin{\triangleright}x)
 =\sum \rho_r(k_{(1)}xS(k_{(2)}))
 =\sum \varepsilon(k_{(2)})\rho_r(k_{(1)}x)
 =k\rho_r(x).
\]
The middle equality holds because
$S(k_{(2)})-\varepsilon(k_{(2)})1\in K^+$ and right multiplication by
this element gives a class in $HK^+$.

\begin{theorem}
\label{thm:symbol-obstruction}
Let $V=\kHom(M,N)$ with its internal-Hom $H$-action.  Evaluation at the
class of $1$ gives natural isomorphisms
\begin{equation}\label{eq:represent-relative}
 \operatorname{Hom}_H(A_r,V)
 \simeq \cD^r_{H\mid K}(M,N).
\end{equation}
For every $r\geq1$, applying $\operatorname{Hom}_H(-,V)$ to
\eqref{eq:relative-layer-sequence} induces an exact sequence
\begin{equation}\label{eq:symbol-obstruction}
 0\longrightarrow
 \gr^r\cD_{H\mid K}(M,N)
 \xrightarrow{\;\bar\sigma_r\;}
 \operatorname{Hom}_H(C_r,V)
 \xrightarrow{\;\partial_r\;}
 \Ext_H^1(A_{r-1},V).
\end{equation}
Under the internal-Hom adjunction, the first arrow becomes the
principal-symbol map of Theorem~\ref{thm:faithful-symbol} after
precomposition with $\rho_r$.  More precisely, the following diagram
commutes:
\[
\begin{CD}
 \gr^r\cD_{H\mid K}(M,N) @>{\bar\sigma_r}>>
       \operatorname{Hom}_H(C_r,V)\\
 @V{\sigma_r}VV @VV{\rho_r^*}V\\
 \operatorname{Hom}_K(Q_r\otimes M,N) @<{\operatorname{adj}}<<
       \operatorname{Hom}_K(Q_r,V^J).
\end{CD}
\]
Here $\rho_r^*(\theta)=\theta\rho_r$, and an $H$-linear map
$C_r\to V$ takes values in $V^J$ because $JC_r=0$.  By a formal
degree-$r$ symbol we mean an element
$\theta\in\operatorname{Hom}_H(C_r,V)$.  Such a symbol is the principal
symbol of a relative higher intertwiner if and only if
$\partial_r(\theta)=0$.  More explicitly, $\partial_r(\theta)$ is the
class of the pushout extension
\[
 0\longrightarrow V\longrightarrow
 (V\oplus A_r)/\{(\theta(c),-c):c\in C_r\}
 \longrightarrow A_{r-1}\longrightarrow0.
\]

In particular, the map
$\bar\sigma_r:\gr^r\cD_{H\mid K}(M,N)\to\operatorname{Hom}_H(C_r,V)$
is surjective if $\Ext_H^1(A_{r-1},V)=0$, for example when $V$ is
injective.  It is also surjective whenever
\eqref{eq:relative-layer-sequence} splits as a sequence of $H$-modules.
These conditions do not in general make $\sigma_r$ surjective onto the
larger space $\operatorname{Hom}_K(Q_r\otimes M,N)$.  Under the
internal-Hom adjunction, its image is exactly
\[
 \operatorname{adj}\!\left(
 \rho_r^*\bigl(\ker\partial_r\bigr)\right)
 \subseteq\operatorname{Hom}_K(Q_r\otimes M,N).
\]
\end{theorem}

\begin{proof}
The sum $HK^++J^{r+1}$ is a left ideal.  An $H$-linear map
$A_r\to V$ is determined by the image $v$ of the class of $1$, and this
image satisfies
\[
 K^+v=0,\qquad J^{r+1}v=0.
\]
Conversely, every $v$ satisfying these equations defines an $H$-linear
map $A_r\to V$.  For the internal-Hom action, $K^+v=0$ is equivalent
to $v\in\operatorname{Hom}_K(M,N)$.  This proves
\eqref{eq:represent-relative}.

Apply the contravariant functor $\operatorname{Hom}_H(-,V)$ to
\eqref{eq:relative-layer-sequence}.  The beginning of the resulting
long exact sequence is
\[
 0\longrightarrow\operatorname{Hom}_H(A_{r-1},V)
 \longrightarrow\operatorname{Hom}_H(A_r,V)
 \longrightarrow\operatorname{Hom}_H(C_r,V)
 \xrightarrow{\partial_r}\Ext_H^1(A_{r-1},V).
\]
Under \eqref{eq:represent-relative}, the first two terms are
$\cD^{r-1}_{H\mid K}(M,N)$ and
$\cD^r_{H\mid K}(M,N)$.  The restriction of the map corresponding to
$f$ sends the class of $x\in J^r$ to $x\cdot f$.  Its composite with
$\rho_r$ is the curried form of $\sigma_r(f)$.  Taking the quotient by
the first term gives \eqref{eq:symbol-obstruction}.

Exactness at $\operatorname{Hom}_H(C_r,V)$ identifies the image with
$\ker\partial_r$.  To describe $\partial_r$, push out the inclusion
$C_r\to A_r$ along $\theta:C_r\to V$.  The resulting middle term is
\[
 (V\oplus A_r)/\{(\theta(c),-c):c\in C_r\},
\]
and its quotient by $V$ is $A_{r-1}$.  This is the extension representing
$\partial_r(\theta)$.

If $\Ext_H^1(A_{r-1},V)=0$, then $\partial_r$ is the zero map, and
exactness makes $\bar\sigma_r$ surjective onto
$\operatorname{Hom}_H(C_r,V)$.  If
\eqref{eq:relative-layer-sequence} splits, every $H$-linear map
$C_r\to V$ extends to $A_r$ by setting it equal to zero on a chosen
copy of $A_{r-1}$.  Its restriction is the prescribed symbol, so the
map $\bar\sigma_r$ is again surjective.  The formula for the image of
$\sigma_r$ follows from the commutative diagram and exactness at
$\operatorname{Hom}_H(C_r,V)$.
\end{proof}

We identify the algebra governing these symbols.  If
$R$ is a Hopf algebra in the braided category
${}_K^K\mathcal{YD}$, its Radford biproduct, or bosonization, is the
ordinary Hopf algebra $R\#K$ \cite{Radford,Majid}.  It has a split Hopf
projection $R\#K\to K$, $r\#k\mapsto\varepsilon_R(r)k$.

\begin{theorem}\label{thm:bosonization}
Let $H=R\#K$ and let $R^+=\ker\varepsilon_R$.  For the kernel $J$ of the
canonical projection $H\to K$ one has, for every $n\geq1$,
\begin{equation}\label{eq:biproduct-powers}
 J^n=(R^+)^n\#K.
\end{equation}
Therefore
\[
 \cD^r_{H\mid K}(M,N)
 =\{f\in\operatorname{Hom}_K(M,N):
       ((R^+)^{r+1}\#1)\cdot f=0\},
\]
and there is a natural graded Hopf algebra isomorphism
\[
 \gr_JH\simeq (\gr_{R^+}R)\#K.
\]
Under this isomorphism, Theorem~\ref{thm:faithful-symbol} identifies the
principal-symbol calculus with convolution by the braided graded
coproduct of $\gr_{R^+}R$ followed by bosonization.
\end{theorem}

\begin{proof}
Throughout the statement, $(R^+)^n$ denotes the $n$th ideal power formed
with the multiplication of the braided algebra object $R$.  As a vector
space it is spanned by products $r_1\cdots r_n$ with $r_i\in R^+$; the
braiding enters the Hopf structure, not the meaning of this product.
As a vector space, the kernel of the projection is $R^+\otimes K$.
The multiplication rule
\[
 (r\#k)(s\#\ell)=r(k_{(1)}\cdot s)\#k_{(2)}\ell
\]
and the $K$-stability of $R^+$ show inductively that
$J^n\subseteq (R^+)^n\#K$.  Conversely, a typical product
$r_1\cdots r_n\#k$, with all $r_i\in R^+$, is
$(r_1\#1)\cdots(r_n\#1)(1\#k)$ and belongs to $J^n$.  This proves
\eqref{eq:biproduct-powers}.  Since right multiplication by $1\#K$
does not change the annihilation condition on a $K$-linear map, the
formula for $\cD^r$ follows.  Taking successive quotients in
\eqref{eq:biproduct-powers} gives the graded algebra isomorphism.  We
describe the coproduct to verify the Hopf assertion.  Write the braided
coproduct of $r\in R$ as $\sum r^{(1)}\otimes r^{(2)}$ and the left
$K$-coaction as $s\mapsto\sum s_{(-1)}\otimes s_{(0)}$.  The biproduct
coproduct is
\[
 \Delta(r\#k)=\sum
 \bigl(r^{(1)}\#(r^{(2)})_{(-1)}k_{(1)}\bigr)
 \otimes\bigl((r^{(2)})_{(0)}\#k_{(2)}\bigr).
\]
The $K$-action and coaction preserve every power of $R^+$.  Moreover,
\[
 \Delta_R((R^+)^n)
 \subseteq\sum_{a+b\geq n}(R^+)^a\otimes(R^+)^b,
\]
by the multiplication argument of Lemma~\ref{lem:filtered-coproduct},
carried out in the braided category.  Consequently the component of total
degree $n$ in the biproduct formula is exactly the bosonization of the
graded braided coproduct on $\gr_{R^+}R$.  This proves the graded Hopf
algebra isomorphism.  The final assertion is
equation~\eqref{eq:symbol-composition}.
\end{proof}

\begin{corollary}
\label{cor:bosonization-symbol-image}
In the setting of Theorem~\ref{thm:bosonization}, put
\[
 Q_r^R=(R^+)^r/(R^+)^{r+1}
 \quad\text{and}\quad
 V^J=\{v\in V:Jv=0\}.
\]
Then $C_r$ is the inflation along $H\to K$ of the $K$-module $Q_r^R$,
and \eqref{eq:symbol-obstruction} becomes
\begin{equation}\label{eq:bosonization-symbol-obstruction}
 0\longrightarrow
 \gr^r\cD_{H\mid K}(M,N)
 \longrightarrow
 \operatorname{Hom}_K(Q_r^R,V^J)
 \xrightarrow{\;\partial_r\;}
 \Ext_H^1(A_{r-1},V).
\end{equation}
Thus the image of the relative principal-symbol map is the kernel of an
explicit Ext obstruction inside the $K$-equivariant braided symbol
space.
\end{corollary}

\begin{proof}
As a right $K$-module, $H=R\#K$ is $R\otimes K$, while
$J^r=(R^+)^r\otimes K$.  Quotienting the latter by $J^{r+1}$ and by
right multiplication with $K^+$ identifies $C_r$ with $Q_r^R$ as a
vector space.  Left multiplication by $R^+$ raises augmentation degree,
so $J$ acts trivially on $C_r$.  Left multiplication by $K$ induces
the given Yetter--Drinfeld $K$-action on $Q_r^R$.  Hence an $H$-linear
map $C_r\to V$ is exactly a $K$-linear map $Q_r^R\to V$ whose image is
annihilated by $J$.  Therefore
\[
 \operatorname{Hom}_H(C_r,V)
 \simeq\operatorname{Hom}_K(Q_r^R,V^J).
\]
Substitution in Theorem~\ref{thm:symbol-obstruction} proves the result.
\end{proof}

\begin{corollary}\label{cor:nichols-symbol}
Suppose $R=\mathcal B(W)$ is a Nichols algebra generated in degree one.
Then relative order for $\mathcal B(W)\#K$ is exactly braided
augmentation order, and the graded Hopf algebra governing formal symbols
retains the defining relations of $\mathcal B(W)$.  Assume in addition
that $k$ has characteristic zero and contains the parameters
$q_i^{\pm1}$, the Cartan matrix is of finite type, and $q$ is not a root
of unity.  Let
\[
 \Gamma=\langle K_1,\ldots,K_\ell\rangle\simeq\mathbb Z^\ell,
 \qquad K=k\Gamma.
\]
Use the coproduct convention
\[
 \Delta(E_i)=E_i\otimes K_i+1\otimes E_i
\]
for the Drinfeld--Jimbo positive Borel, and set
\[
 R=U_q(\mathfrak b^+)^{\operatorname{co}\pi_+}
  =\{x\in U_q(\mathfrak b^+):(\id\otimes\pi_+)\Delta(x)=x\otimes1\}.
\]
Put
$e_i=E_iK_i^{-1}$.  The right coinvariant algebra for the projection to
$K$ is $R$, and $R$ is the Nichols algebra $\mathcal B(W)$ on the diagonal
Yetter--Drinfeld module $W$ with basis $e_1,\ldots,e_\ell$, action and
coaction
\[
 K_i\cdot e_j=q_i^{a_{ij}}e_j,
 \qquad \delta(e_i)=K_i^{-1}\otimes e_i,
\]
and braiding
\begin{equation}\label{eq:inverse-cartan-braiding}
 c(e_i\otimes e_j)=q_i^{-a_{ij}}e_j\otimes e_i.
\end{equation}
Its bosonization by $k\Gamma$ maps isomorphically onto the positive Borel
by
\[
 e_i\#1\longmapsto E_iK_i^{-1},
 \qquad 1\#K_j\longmapsto K_j;
\]
equivalently, $e_i\#K_i$ maps to $E_i$.  Thus the quantum Serre relations
with inverse parameters constrain the mixed symbols for connected Dynkin
diagrams, and the principal-symbol image is selected by the obstruction
map in Corollary~\ref{cor:bosonization-symbol-image}.
\end{corollary}

\begin{proof}
A Nichols algebra is connected graded and generated by $W=R_1$, so
$(R^+)^n=\bigoplus_{d\geq n}R_d$.  Apply
Theorem~\ref{thm:bosonization}.  We verify the convention in the Borel
specialization.  Since
\[
 \Delta(e_i)=e_i\otimes1+K_i^{-1}\otimes e_i,
\]
one has $(\id\otimes\pi_+)\Delta(e_i)=e_i\otimes1$, so $e_i$ is a right
coinvariant.  The left coaction obtained from the Radford decomposition is
$e_i\mapsto K_i^{-1}\otimes e_i$.  Together with the adjoint Cartan action,
this gives \eqref{eq:inverse-cartan-braiding}.  The multiplication map from
the bosonization to the Borel sends $e_i\#K_i$ to $E_i$, and therefore
contains all the standard Borel generators.

In finite Cartan type, the defining Nichols relations are the quantum
Serre relations for the inverse Cartan parameters.  Replacing every
$q_i$ by $q_i^{-1}$ preserves those relations, after the corresponding
change of braided convention.  Hence the multiplication map is an
isomorphism by the PBW theorem.  See \cite{Rosso},
\cite[Section~1]{AndruskiewitschSchneider}, and
\cite[Chapter~8]{Jantzen}.  The cited PBW construction is usually stated
over a rational-function field.  Under the standing hypotheses it
specializes to $k$ by base change, since the relevant quantum integers do
not vanish.  At roots of unity further power relations may
occur, and no root-of-unity identification is asserted here.
\end{proof}

\subsection{Positive and negative quantum Borels}

For the remainder of the quantum-Borel calculations, $k$ has
characteristic zero, the Cartan matrix is of finite type, and $q$ is not a
root of unity.  We use the conventions fixed in Section~4.3.

Let $U_q(\mathfrak b^+)$ be the positive Borel subalgebra generated by
$K_i^{\pm1}$ and $E_i$.  Let $U_q(\mathfrak h)$ be the group algebra
generated by the $K_i^{\pm1}$.  There is a split Hopf projection
\[
 \pi_+:U_q(\mathfrak b^+)\longrightarrow U_q(\mathfrak h),
 \qquad E_i\longmapsto0,\quad K_i\longmapsto K_i.
\]
Its kernel $J_+$ is the ideal generated by the $E_i$.  The negative Borel
has the analogous projection $\pi_-$ with kernel $J_-$ generated by the
$F_i$.

For a positive root $\beta=\sum_i n_i\alpha_i$, write
$\operatorname{ht}(\beta)=\sum_i n_i$.  Give a PBW monomial in positive
root vectors the sum of the heights of its factors.

\begin{proposition}\label{prop:root-height}
With the hypotheses above, let $U_q(\mathfrak b^+)_{\rm ht}$ denote the
same Hopf algebra equipped with the root-height grading.  Then $J_+^r$ is
spanned by the
PBW monomials of total positive-root height at least $r$, with arbitrary
Cartan factors.  Similarly, $J_-^r$ is spanned by negative-root PBW
monomials of total height at least $r$.
Moreover, root height gives an $\mathbb N$-grading
\[
 U_q(\mathfrak b^+)=\bigoplus_{d\geq0}B_d^+,
 \qquad J_+^r=\bigoplus_{d\geq r}B_d^+,
\]
and the canonical map induced by homogeneous representatives is an
isomorphism of graded Hopf algebras
\begin{equation}\label{eq:graded-borel}
 \gr_{J_+}U_q(\mathfrak b^+)\simeq U_q(\mathfrak b^+)_{\rm ht}.
\end{equation}
The corresponding assertions hold for the negative Borel.
\end{proposition}

\begin{proof}
Assign degree zero to the $K_i^{\pm1}$ and degree one to every $E_i$.
The Cartan-root commutation relations and the quantum Serre relations are
homogeneous, so this assignment defines an $\mathbb N$-grading.  We use the
PBW basis obtained from a reduced expression of the longest Weyl-group
element, with its associated ordered quantum root vectors; see
\cite[Chapter~8]{Jantzen}.  The PBW theorem identifies the degree-$d$
component $B_d^+$ with the span of the PBW monomials whose total root
height is $d$.  The projection $\pi_+$ is
the projection onto $B_0^+$, and hence
$J_+=\bigoplus_{d\geq1}B_d^+$.  A product of $r$ elements of $J_+$ has
degree at least $r$, so
$J_+^r\subseteq\bigoplus_{d\geq r}B_d^+$.

For the reverse inclusion, a positive root vector of height $d$ is a
linear combination of words of length $d$ in the simple generators.
After commuting Cartan factors past root vectors, every PBW monomial of
total height $d$ therefore belongs to $J_+^d$, and thus to $J_+^r$ when
$d\geq r$.  This proves the equality.  The map which sends the class of
$x\in B_d^+$ modulo $J_+^{d+1}$ to $x$ is consequently a graded vector
space isomorphism.  Multiplication is preserved because degrees add.
The coproduct is also homogeneous for total root height, as follows on
the generators from
$\Delta(E_i)=E_i\otimes K_i+1\otimes E_i$; hence the map is a Hopf
algebra isomorphism.  Assigning degree one to each $F_i$ and degree zero
to each $K_i^{\pm1}$ gives the corresponding statement for the negative
Borel.
\end{proof}

\begin{corollary}\label{cor:root-symbol}
For any $U_q(\mathfrak b^+)$-modules $M,N$ and every $r\geq0$, there is a
natural injection
\[
 \gr^r\cD_{U_q(\mathfrak b^+)\mid U_q(\mathfrak h)}(M,N)
 \lhook\joinrel\longrightarrow
 \operatorname{Hom}_{U_q(\mathfrak h)}(B_r^+\otimes M,N),
\]
where $B_r^+$ is the root-height-$r$ component of the positive quantum
Borel.  For $r\geq1$, its image is described as follows.  Put
\[
 H_+=U_q(\mathfrak b^+),\qquad K=U_q(\mathfrak h),\qquad
 V=\kHom(M,N),
\]
and let $K^+=\ker(\varepsilon|_K)$.  Under the identification
$B_r^+\simeq J_+^r/J_+^{r+1}$ of
Proposition~\ref{prop:root-height}, define
\[
 C_r=\frac{H_+K^++J_+^r}{H_+K^++J_+^{r+1}},\qquad
 A_{r-1}=\frac{H_+}{H_+K^++J_+^r},
\]
and let $\rho_r:B_r^+\to C_r$ be the natural quotient map.  Precomposition
gives
\[
 \rho_r^*: \operatorname{Hom}_{H_+}(C_r,V)
 \longrightarrow \operatorname{Hom}_K(B_r^+,V^{J_+}),
 \qquad \theta\longmapsto\theta\rho_r,
\]
while internal-Hom adjunction sends $\varphi:B_r^+\to V^{J_+}$ to
the map $\operatorname{adj}(\varphi)$ given by
\[
 \operatorname{adj}(\varphi)(x\otimes m)=\varphi(x)(m).
\]
If
\[
 \partial_r:\operatorname{Hom}_{H_+}(C_r,V)
 \longrightarrow \Ext_{H_+}^1(A_{r-1},V)
\]
is the connecting map associated with
$0\to C_r\to A_r\to A_{r-1}\to0$, then the image of the principal-symbol
injection is exactly
\[
 \operatorname{adj}\!\left(\rho_r^*(\ker\partial_r)\right).
\]
Thus a symbol is liftable precisely when the corresponding map
$C_r\to V$ has zero connecting class.  Composition of liftable symbols is
given by the convolution formula
\eqref{eq:symbol-composition}, using the $(p,q)$ component of the
root-height graded coproduct.  The analogous statement holds for the
negative Borel.
\end{corollary}

\begin{proof}
Proposition~\ref{prop:root-height} identifies
$J_+^r/J_+^{r+1}$ with $B_r^+$.  The injection and its convolution law are
therefore Theorem~\ref{thm:faithful-symbol}.  The description of its image
is Theorem~\ref{thm:symbol-obstruction}.  The negative case is identical.
\end{proof}

\begin{corollary}\label{cor:borel-orders}
For $U_q(\mathfrak g)$-modules $M,N$, a Cartan-equivariant map $f$ lies
in
\[
 \cD^r_{U_q(\mathfrak b^+)\mid U_q(\mathfrak h)}(M,N)
\]
if and only if
\[
 (E_{i_0}\cdot)\,(E_{i_1}\cdot)\cdots(E_{i_r}\cdot)f=0
\]
for all simple roots $i_0,\ldots,i_r$.  The analogous statement holds
for the negative Borel with the $F_i$.
\end{corollary}

\begin{proof}
The ideal $J_+^{r+1}$ is spanned by words of length at least $r+1$ in the
$E_i$, with Cartan coefficients.  Cartan equivariance of $f$ and
Proposition~\ref{prop:root-height} reduce the annihilation condition to
the displayed words.  The negative case is the same.
\end{proof}

Define a bifiltration on Cartan-equivariant maps by
\[
 \cD^{p,q}_q(M,N)=
 \cD^p_{U_q(\mathfrak b^-)\mid U_q(\mathfrak h)}(M,N)
 \cap
 \cD^q_{U_q(\mathfrak b^+)\mid U_q(\mathfrak h)}(M,N).
\]

\begin{theorem}\label{thm:bifiltration}
For $U_q(\mathfrak g)$-modules $M,N,P$,
\[
 \cD^{p,q}_q(N,P)\circ\cD^{p',q'}_q(M,N)
 \subseteq
 \cD^{p+p',q+q'}_q(M,P).
\]
Moreover,
\[
 \cD^{0,0}_q(M,N)=\operatorname{Hom}_{U_q(\mathfrak g)}(M,N).
\]
If for each Cartan-equivariant $f:M\to N$ there are $p,q\geq0$ such that
$J_-^{p+1}\cdot f=0$ and $J_+^{q+1}\cdot f=0$, then the bifiltration is
exhaustive on $\operatorname{Hom}_{U_q(\mathfrak h)}(M,N)$.
\end{theorem}

\begin{proof}
Apply Theorem~\ref{thm:relative-composition} separately to the positive
and negative split Hopf pairs.  A map in bidegree $(0,0)$ intertwines the
Cartan generators and all $E_i,F_i$, hence is
$U_q(\mathfrak g)$-linear.  Conversely, a $U_q(\mathfrak g)$-linear map
intertwines the Cartan generators and every $E_i$ and $F_i$, so it has
relative order zero for both Borels.

Under the additional hypothesis in the last sentence, every
Cartan-equivariant map belongs to one of the displayed intersections,
which is precisely exhaustivity.  This condition alone does not imply
noncollapse.  Noncollapse is equivalent to the existence of a nonzero
quotient in positive bidegree.  The modules in the following subsections
give such quotients for the one-sided Borel filtrations.
\end{proof}

\subsection{Strictness in rank one}

We show that the relative filtration has nonzero symbols in every degree.
Let $B=U_q(\mathfrak b^+_{\mathfrak{sl}_2})$, with
\[
 KEK^{-1}=q^2E,\qquad
 \Delta(E)=E\otimes K+1\otimes E.
\]
Let $V=k[t]$ be the $B$-module
\[
 Kt^n=q^{2n}t^n,\qquad Et^n=t^{n+1}.
\]
For a polynomial $P(z)\in k[z]$, define the Cartan-equivariant diagonal
map
\[
 T_P(t^n)=P(q^{2n})t^n.
\]

\begin{proposition}\label{prop:strict-rank-one}
Assume that $q$ is not a root of unity.  If $P$ is nonzero and has degree
$d$, then
\[
 T_P\in\cD^d_{B\mid k[K^{\pm1}]}(V,V)
 \setminus
 \cD^{d-1}_{B\mid k[K^{\pm1}]}(V,V).
\]
In particular, every nonnegative relative order occurs.
\end{proposition}

\begin{proof}
More generally, suppose $T(t^n)=a_n t^{n+s}$.  Using
$S(E)=-EK^{-1}$ and the internal-Hom action gives
\[
 (E\cdot T)(t^n)
 =q^{-2n}(a_n-a_{n+1})t^{n+s+1}.
\]
For $a_n=P(q^{2n})$, put $z=q^{2n}$.  The action of $E$ replaces $P$ by
\[
 \partial_qP(z)=z^{-1}\bigl(P(z)-P(q^2z)\bigr).
\]
More precisely, induction on $j$ gives
\[
 (E^j\cdot T_P)(t^n)=(\partial_q^jP)(q^{2n})t^{n+j}.
\]
For a monomial one has
\[
 \partial_q(z^m)=(1-q^{2m})z^{m-1}.
\]
The coefficient is nonzero for $m>0$ because $q$ is not a root of unity.
It follows that $\partial_q^{d+1}P=0$ and
$\partial_q^dP\ne0$ when $P$ has degree $d$.  Hence
$E^{d+1}\cdot T_P=0$ while $E^d\cdot T_P\ne0$.  Since the kernel
ideal of $B\to k[K^{\pm1}]$ is generated by $E$, Corollary
\ref{cor:borel-orders} proves the assertion.
\end{proof}

\begin{theorem}\label{thm:rank-one-classification}
Assume that $q$ is not a root of unity.  For every $r\ge0$,
\[
 \cD^r_{B\mid k[K^{\pm1}]}(V,V)
 =\{T_P:P\in k[z],\ \deg P\le r\}.
\]
Consequently,
\[
 \cD^\infty_{B\mid k[K^{\pm1}]}(V,V)\cong k[z]
\]
as filtered algebras, where composition corresponds to multiplication of
polynomials and relative order is polynomial degree.
\end{theorem}

\begin{proof}
The $K$-weights $q^{2n}$ of $V$ are pairwise distinct.  Every
$K$-linear endomorphism is therefore diagonal, say
$T(t^n)=a_nt^n$.  Put
\[
 (\delta a)_n=q^{-2n}(a_n-a_{n+1}).
\]
The calculation in Proposition~\ref{prop:strict-rank-one} shows that
$E^s\cdot T=0$ if and only if $\delta^sa=0$.  We prove by induction on
$r$ that $\delta^{r+1}a=0$ precisely when
$a_n=P(q^{2n})$ for a polynomial $P$ of degree at most $r$.

For $r=0$, the equation $\delta a=0$ says that the sequence is constant.
Suppose the assertion is known in degree $r-1$ and
$\delta^{r+1}a=0$.  Then
\[
 (\delta a)_n=Q(q^{2n})
\]
for some $Q$ of degree at most $r-1$.  Equivalently,
\[
 a_n-a_{n+1}=q^{2n}Q(q^{2n}).
\]
Write $Q(z)=\sum_{j=0}^{r-1}c_jz^j$.  Since
$1-q^{2(j+1)}\ne0$, the polynomial
\[
 R(z)=\sum_{j=0}^{r-1}\frac{c_j}{1-q^{2(j+1)}}z^{j+1}
\]
has degree at most $r$, has zero constant term, and satisfies
\[
 R(q^{2n})-R(q^{2n+2})=q^{2n}Q(q^{2n}).
\]
Hence $a_n-R(q^{2n})$ is constant.  This proves the forward implication.
For the reverse implication, repeated application of $\delta$ lowers
polynomial degree until the sequence becomes zero.  Since
$T_P\circ T_Q=T_{PQ}$, the identification preserves
the filtrations and products.
\end{proof}

\subsection{Simple-root strings in arbitrary rank}

The rank-one calculation is not confined to algebras of rank one.  It
occurs along each simple-root direction of an arbitrary Drinfeld--Jimbo
Borel and supplies a uniform strictness result for the root-height
filtration.  Write $(a_{ji})$ for the Cartan matrix and $q_j=q^{d_j}$,
with $d_ja_{ji}=d_ia_{ij}$.

Fix a simple root $i$, and let $V_i=k[t]$.  Define an action of
$U_q(\mathfrak b^+)$ by
\begin{equation}\label{eq:simple-root-string}
 K_jt^n=q_j^{a_{ji}n}t^n,\qquad
 E_it^n=t^{n+1},\qquad E_jt^n=0\quad(j\ne i).
\end{equation}
The Cartan relations hold because
$K_jE_iK_j^{-1}=q_j^{a_{ji}}E_i$, and every quantum Serre relation is
zero: it either involves an $E_j$ with $j\ne i$, or is vacuous.  Hence
\eqref{eq:simple-root-string} is a Borel module.

\begin{theorem}
\label{thm:simple-root-classification}
Assume that $q_i^2$ is not a root of unity.  For $P\in k[z]$ define
\[
 T_P(t^n)=P(q_i^{2n})t^n.
\]
Then
\[
 \cD^r_{U_q(\mathfrak b^+)\mid U_q(\mathfrak h)}(V_i,V_i)
   =\{T_P:\deg P\le r\}.
\]
Consequently its union over $r$ is $k[z]$ as a filtered algebra.  In
particular, for every finite rank, every simple root, and every $r\geq0$,
the root-height filtration has a nonzero symbol of order $r$, and the
composition bound in Theorem~\ref{thm:relative-composition} is sharp.
\end{theorem}

\begin{proof}
The $K_i$-weights of $V_i$ are the pairwise distinct scalars
$q_i^{2n}$.  Thus every $U_q(\mathfrak h)$-linear endomorphism is
diagonal.  On such an endomorphism the internal-Hom action of $E_j$ is
zero for $j\ne i$, since $E_j$ acts as zero on both source and target.
If $T(t^n)=a_nt^n$, the coproduct and antipode formulas for $E_i$ give
\[
 (E_i\cdot T)(t^n)=q_i^{-2n}(a_n-a_{n+1})t^{n+1}.
\]
Thus the $E_i$-defect is the one-variable quantum difference used in
Theorem~\ref{thm:rank-one-classification}.  Its kernel after $r+1$
iterations consists exactly of the sequences
$a_n=P(q_i^{2n})$ with $\deg P\le r$.  The other simple-root defects
vanish, so Corollary~\ref{cor:borel-orders} gives the displayed equality.
Taking $P=z^r$ proves strictness, while
$T_PT_Q=T_{PQ}$ proves sharpness under composition.
\end{proof}

Theorem~\ref{thm:simple-root-classification} classifies one root direction.
A tensor product of these strings does
not automatically define a module for a connected quantum Borel:
adjacent simple-root operators are constrained by the quantum Serre
relations, and their iterated defects need not be independent commuting
finite differences.  Thus PBW root height identifies the symbol spaces
and proves strictness, but does not by itself imply a multivariable polynomial
classification.  Such a classification requires additional information
about the module and its mixed-root defect operators.

For mutually orthogonal roots, the quantum Serre relations do not couple
the defect operators.  We first record the elementary
finite-difference fact needed for a complete classification.  If
$Q_1,\ldots,Q_m\in k^\times$ are not roots of unity, put
\[
 (\delta_j a)(\mathbf n)=Q_j^{-n_j}
 \bigl(a(\mathbf n)-a(\mathbf n+\mathbf e_j)\bigr)
\]
for a function $a:\mathbb N^m\to k$.

\begin{lemma}
\label{lem:multivariable-differences}
The following conditions on $a:\mathbb N^m\to k$ are equivalent:
\begin{enumerate}
\item $\delta_{i_1}\cdots\delta_{i_{r+1}}a=0$ for every sequence
$i_1,\ldots,i_{r+1}$;
\item there is a polynomial $P\in k[z_1,\ldots,z_m]$ of total degree at
most $r$ such that
\[
 a(\mathbf n)=P(Q_1^{n_1},\ldots,Q_m^{n_m}).
\]
\end{enumerate}
The polynomial $P$ is unique.
\end{lemma}

\begin{proof}
On functions of $z_j=Q_j^{n_j}$, the operator $\delta_j$ is
\[
 z_j^{-1}\bigl(P(\mathbf z)-P(z_1,\ldots,Q_jz_j,\ldots,z_m)\bigr).
\]
It lowers the degree in $z_j$ by one and has a nonzero scalar on every
positive power of $z_j$.  Hence (2) implies (1).

For the converse, argue by induction on $m$.  The case $m=1$ is the
antidifference argument in Theorem~\ref{thm:rank-one-classification}.
Only division by the nonzero scalars $1-Q_m^j$ occurs in that argument,
so the same calculation works for sequences with values in a vector
space.  Apply it in the last variable, using the space of functions on
$\mathbb N^{m-1}$ as the coefficient space.  The equation
$\delta_m^{r+1}a=0$ gives
\[
 a(\mathbf n)=\sum_{d=0}^r b_d(n_1,\ldots,n_{m-1})Q_m^{dn_m}.
\]
Fix $d$ and let $W$ be any word of length $r-d+1$ in
$\delta_1,\ldots,\delta_{m-1}$.  The operators in different variables
commute.  The hypothesis therefore gives $W\delta_m^d a=0$.  In the
displayed expansion, $\delta_m^d$ kills the terms with exponent smaller
than $d$.  The term with exponent $d$ becomes a nonzero scalar multiple
of $b_d$, while every term with larger exponent retains positive degree
in $Q_m^{n_m}$.  The constant coefficient of $W\delta_m^d a$ is thus a
nonzero scalar multiple of $Wb_d$.  Since a polynomial that vanishes at
all $Q_m^{n_m}$ is zero, $Wb_d=0$.  Thus $b_d$
satisfies condition (1) in $m-1$ variables with order $r-d$.  Induction
shows that $b_d$ is a polynomial of total degree at most $r-d$ in
$Q_1^{n_1},\ldots,Q_{m-1}^{n_{m-1}}$.  This proves (2).  Uniqueness
follows successively from the infinitude of each set
$\{Q_j^n:n\geq0\}$.
\end{proof}

Let $S=\{i_1,\ldots,i_m\}$ be a set of pairwise orthogonal simple roots,
so $a_{i_s i_t}=0$ for $s\ne t$.  Let $V_S$ have basis
$v_{\mathbf n}$, $\mathbf n\in\mathbb N^m$, and define
\[
 K_jv_{\mathbf n}
 =\left(\prod_{s=1}^m q_j^{a_{j i_s}n_s}\right)v_{\mathbf n},
 \qquad E_{i_s}v_{\mathbf n}=v_{\mathbf n+\mathbf e_s},
 \qquad E_jv_{\mathbf n}=0\quad(j\notin S).
\]
The quantum Serre relations hold: the active generators commute because
the roots in $S$ are orthogonal, and every remaining Serre relation has
an inactive generator.

\begin{theorem}
\label{thm:orthogonal-root-classification}
Assume that $q$ is not a root of unity and that the finite-type Cartan
matrix is nonsingular.  For $P\in k[z_1,\ldots,z_m]$, define
\[
 T_P(v_{\mathbf n})
 =P(q_{i_1}^{2n_1},\ldots,q_{i_m}^{2n_m})v_{\mathbf n}.
\]
Then
\[
 \cD^r_{U_q(\mathfrak b^+)\mid U_q(\mathfrak h)}(V_S,V_S)
 =\{T_P:\deg_{\rm tot}P\leq r\}.
\]
Consequently the union of the relative filtration is the polynomial
algebra $k[z_1,\ldots,z_m]$, filtered by total degree.  Composition is
multiplication, and the order estimate is sharp in every multidegree.
\end{theorem}

\begin{proof}
Nonsingularity of the Cartan matrix implies that the joint Cartan weights
of the $v_{\mathbf n}$ are distinct.  Hence every Cartan-linear
endomorphism is diagonal, say $T(v_{\mathbf n})=a(\mathbf n)v_{\mathbf n}$.
For an active generator $E_{i_s}$, orthogonality gives
\[
 (E_{i_s}\cdot T)(v_{\mathbf n})
 =q_{i_s}^{-2n_s}
 \bigl(a(\mathbf n)-a(\mathbf n+\mathbf e_s)\bigr)
 v_{\mathbf n+\mathbf e_s}.
\]
The inactive generators act trivially on the internal Hom.  By
Corollary~\ref{cor:borel-orders}, $T$ has order at most $r$ exactly when
every word of length $r+1$ in these commuting finite differences
annihilates $a$.  Lemma~\ref{lem:multivariable-differences}, with
$Q_s=q_{i_s}^2$, gives the asserted classification.  The identity
$T_PT_Q=T_{PQ}$ identifies composition with polynomial multiplication,
and total degree is additive on nonzero products over a field.
\end{proof}

\subsection{Connected simply laced types}

We give a complete calculation with all simple-root directions active.
The answer is controlled by directed paths in the Dynkin diagram.  Let
$A=(a_{ij})$ be a simply laced Cartan matrix of finite type and choose an
orientation of its Dynkin diagram.  A finite-type Dynkin diagram is a tree,
so every orientation is acyclic and has a topological ordering.  After
relabeling the vertices by
$1,\ldots,\ell$, we may and do assume that every arrow points from its
smaller to its larger endpoint.  Write $i\sim j$ when the corresponding vertices are
adjacent and put
\[
 \lambda_i(\mathbf n)=\sum_{j<i,\,j\sim i}n_j,
 \qquad \mathbf n=(n_1,\ldots,n_\ell)\in\mathbb N^\ell.
\]
Let $V$ have basis $v_{\mathbf n}$ and define
\begin{equation}\label{eq:oriented-quantum-space-module}
 K_i v_{\mathbf n}=q^{\sum_j a_{ij}n_j}v_{\mathbf n},
 \qquad
 E_i v_{\mathbf n}=q^{\lambda_i(\mathbf n)}
 v_{\mathbf n+\mathbf e_i}.
\end{equation}
\begin{lemma}\label{lem:oriented-module}
The formulas \eqref{eq:oriented-quantum-space-module} define a
$U_q(\mathfrak b^+)$-module.
\end{lemma}

\begin{proof}
The operators $K_i$ are invertible and commute.  Applying
$K_hE_iK_h^{-1}$ to $v_{\mathbf n}$ gives
\[
 K_hE_iK_h^{-1}=q^{a_{hi}}E_i,
\]
which is the Cartan-root relation.  If $i<j$ and $i\sim j$, the term
$n_i$ occurs in $\lambda_j$ but no term involving $n_j$ occurs in
$\lambda_i$.  Consequently $E_iE_j=q^{-1}E_jE_i$.  If $i$ and $j$ are
not adjacent, the two operators commute.

It remains to check both quantum Serre relations along an edge.  Put
$X=E_i$ and $Y=E_j$ for $i<j$.  Since $XY=q^{-1}YX$, one has
\begin{align*}
 X^2Y-(q+q^{-1})XYX+YX^2
 &=\bigl(q^{-2}-(q+q^{-1})q^{-1}+1\bigr)YX^2=0,\\
 Y^2X-(q+q^{-1})YXY+XY^2
 &=\bigl(q^2-(q+q^{-1})q+1\bigr)XY^2=0.
\end{align*}
These are all the Serre relations in simply laced type, so the formulas
define the asserted Borel module.
\end{proof}

Define positive integers recursively by
\begin{equation}\label{eq:path-weights}
 w_i=1+\sum_{j<i,\,j\sim i}w_j.
\end{equation}
Equivalently, $w_i$ is the number of directed paths ending at $i$, with
the length-zero path included.  For
$P\in k[X_1,\ldots,X_\ell]$, put
\[
 T_P(v_{\mathbf n})
 =P(q^{2n_1},\ldots,q^{2n_\ell})v_{\mathbf n}
\]
and let $\operatorname{wtdeg}P$ denote weighted degree for
$\deg X_i=w_i$.

\begin{theorem}
\label{thm:directed-path-classification}
Assume that $q$ is not a root of unity.  For the module
\eqref{eq:oriented-quantum-space-module} and every $r\geq0$,
\begin{equation}\label{eq:directed-path-classification}
 \cD^r_{U_q(\mathfrak b^+)\mid U_q(\mathfrak h)}(V,V)
 =\{T_P:\operatorname{wtdeg}P\leq r\}.
\end{equation}
Consequently,
\[
 \cD^\infty_{U_q(\mathfrak b^+)\mid U_q(\mathfrak h)}(V,V)
 \cong k[X_1,\ldots,X_\ell]
\]
as filtered algebras.  Every order occurs, and exact order is additive
under composition of nonzero elements.
\end{theorem}

\begin{proof}
The joint Cartan weights of the vectors $v_{\mathbf n}$ are distinct.
Indeed, equality of two such weights and the assumption on $q$ give
$A\mathbf n=A\mathbf n'$, while a finite-type Cartan matrix is
nonsingular.  Every Cartan-equivariant endomorphism is therefore diagonal.

We first compute the exact order of a monomial.  For
$\mathbf u\in\mathbb N^\ell$ and $\mathbf A\in\mathbb Z^\ell$, write
\[
 R_{\mathbf u}^{\mathbf A}(v_{\mathbf n})
 =q^{\mathbf A\cdot\mathbf n}v_{\mathbf n+\mathbf u}.
\]
The exponent vector is allowed to lie in $\mathbb Z^\ell$ because
successive defects subtract Cartan rows and can create negative entries.
The expression remains well defined for all such vectors since
$q\in k^\times$.
Let $\boldsymbol\alpha_i=(a_{i1},\ldots,a_{i\ell})$ and let
$\boldsymbol\chi_i=\sum_{j<i,\,j\sim i}\mathbf e_j$.  From
$\Delta(E_i)=E_i\otimes K_i+1\otimes E_i$ and
$S(E_i)=-E_iK_i^{-1}$ one obtains
\begin{equation}\label{eq:general-defect-shift}
 E_i\cdot R_{\mathbf u}^{\mathbf A}
 =\bigl(q^{\lambda_i(\mathbf u)}-q^{A_i}\bigr)
 R_{\mathbf u+\mathbf e_i}^{\mathbf A-
   \boldsymbol\alpha_i+\boldsymbol\chi_i}.
\end{equation}
To see the coefficient in this formula, apply the two coproduct terms to
$v_{\mathbf n}$.  The first contributes the scalar
$q^{\lambda_i(\mathbf u)}$, because $E_i$ acts after the shift by
$\mathbf u$.  The second contributes $-q^{A_i}$, because applying
$R_{\mathbf u}^{\mathbf A}$ after $E_iK_i^{-1}$ increases the exponent
of $q$ by $A_i$.  Both terms have exponent vector
$\mathbf A-\boldsymbol\alpha_i+\boldsymbol\chi_i$, which gives
\eqref{eq:general-defect-shift}.

Start with $T_{\mathbf X^{\mathbf a}}=R_{\mathbf0}^{2\mathbf a}$.
After a nonzero word containing $u_i$ occurrences of $E_i$, the $i$th
component of the exponent vector is
\[
 A_i=2a_i-2u_i+
 \sum_{j<i,\,j\sim i}u_j+
 2\sum_{j>i,\,j\sim i}u_j.
\]
Here an occurrence of $E_i$ lowers the $i$th component by $2$.  An
occurrence of an adjacent $E_j$ raises it by $1$ when $j<i$ and by $2$
when $j>i$; the second contribution in the latter case comes from
$\boldsymbol\chi_j$.
The exponent vector after any prefix of the word depends only on the
numbers of occurrences in that prefix.  The order of the letters matters
through the nonvanishing condition at each successive prefix.
By \eqref{eq:general-defect-shift}, the next $E_i$ action vanishes
precisely when
\begin{equation}\label{eq:capacity}
 u_i=a_i+\sum_{j>i,\,j\sim i}u_j.
\end{equation}
Define $u_i^*(\mathbf a)$, in decreasing order of $i$, by
\begin{equation}\label{eq:maximal-counts}
 u_i^*(\mathbf a)=a_i+
 \sum_{j>i,\,j\sim i}u_j^*(\mathbf a).
\end{equation}
Every word whose successive action on
$T_{\mathbf X^{\mathbf a}}$ remains nonzero has
$u_i\leq u_i^*(\mathbf a)$ for all $i$.  We prove
this by descending induction on $i$.  Before each occurrence of $E_i$ in
a prefix, the current number of $E_i$ letters is strictly smaller than the
right side of \eqref{eq:capacity}; otherwise that occurrence would
annihilate the result of the preceding prefix.  Replacing the current
counts at vertices $j>i$ by their final counts can only increase the right
side.  Once the inductive bounds for those later vertices are known,
\eqref{eq:maximal-counts} bounds the total number of occurrences of $E_i$.
Conversely, apply the generators in decreasing
vertex order, using $E_i$ exactly $u_i^*(\mathbf a)$ times.  At every
stage the inequality preceding \eqref{eq:capacity} is strict, so the
result is nonzero.  Hence
\begin{equation}\label{eq:monomial-path-order}
 \operatorname{ord}(T_{\mathbf X^{\mathbf a}})
 =\sum_i u_i^*(\mathbf a)=\sum_i w_i a_i.
\end{equation}
The last equality follows from \eqref{eq:maximal-counts}: the coefficient
of $a_i$ counts the directed paths ending at $i$, which is the recursion
\eqref{eq:path-weights}.

The same formula holds for a nonzero polynomial without cancellation.
Choose a monomial $\mathbf X^{\mathbf a}$ of maximal weighted degree in
$P$, and apply the word
\[
 E_1^{u_1^*(\mathbf a)}E_2^{u_2^*(\mathbf a)}\cdots
 E_\ell^{u_\ell^*(\mathbf a)},
\]
whose rightmost block acts first.  Repeated use of
\eqref{eq:capacity} shows that this word is nonzero on
$T_{\mathbf X^{\mathbf a'}}$ only if $a_i'\geq a_i$ for every $i$.
Indeed, when the $E_i$ block is applied, all later blocks already have
the fixed sizes $u_j^*(\mathbf a)$, and survival is equivalent to
\[
 u_i^*(\mathbf a)\leq a_i'+
 \sum_{j>i,\,j\sim i}u_j^*(\mathbf a),
\]
which reduces to $a_i\leq a_i'$.  If these inequalities hold and
$\mathbf a'\ne\mathbf a$, positivity of the $w_i$ gives
$\sum_iw_i a_i'>\sum_iw_i a_i$, contrary to the maximal choice of
$\mathbf a$.  Thus no other monomial in the support of $P$ has this
property.  The
chosen monomial therefore cannot cancel, and
\begin{equation}\label{eq:polynomial-path-order}
 \operatorname{ord}(T_P)=\operatorname{wtdeg}P.
\end{equation}

It remains to prove exhaustivity.  Let $T\in\cD^r(V,V)$.  It is diagonal,
say $T(v_{\mathbf n})=c(\mathbf n)v_{\mathbf n}$.  The element
$E_i^{r+1}$ annihilates $T$ for every $i$.  With all coordinates except
$n_i$ fixed, the corresponding $i$-string differs from the rank-one
module only by nonzero constant factors in the actions of $K_i$ and
$E_i$.  The antidifference argument of
Theorem~\ref{thm:rank-one-classification} shows that
$c(\mathbf n)$ is a polynomial of degree at most $r$ in $q^{2n_i}$,
separately for each $i$.

We spell out the interpolation step.  For fixed
$n_2,\ldots,n_\ell$, the values at $n_1=0,\ldots,r$ determine the
unique polynomial of degree at most $r$ in $X_1=q^{2n_1}$.  Its
coefficients are fixed linear combinations of those $r+1$ values and
are therefore separately polynomial functions of
$q^{2n_2},\ldots,q^{2n_\ell}$.  Repeating the same argument for
$n_2,\ldots,n_\ell$ produces one polynomial
$P\in k[X_1,\ldots,X_\ell]$, of degree at most $r$ in each variable,
such that $c(\mathbf n)=P(q^{2n_1},\ldots,q^{2n_\ell})$ for every
$\mathbf n$.  This first step proves polynomiality and separate degree
bounds, but it does not yet give the weighted total-degree bound.  The
exact-order formula \eqref{eq:polynomial-path-order} supplies the second
step: since $T=T_P$ has order at most $r$, it forces
$\operatorname{wtdeg}P\leq r$, proving
\eqref{eq:directed-path-classification}.  Since
$T_PT_Q=T_{PQ}$ and weighted degree is additive on nonzero products,
exact order is additive under composition.
\end{proof}

\begin{corollary}
\label{cor:path-hilbert-series}
For the module in Theorem~\ref{thm:directed-path-classification},
\[
 \gr\cD_{U_q(\mathfrak b^+)\mid U_q(\mathfrak h)}(V,V)
 \cong k[\xi_1,\ldots,\xi_\ell],
 \qquad \deg\xi_i=w_i.
\]
Its Hilbert series is
\begin{equation}\label{eq:path-hilbert-series}
 \sum_{r\geq0}\dim_k\gr^r\cD(V,V)z^r
 =\prod_{i=1}^\ell\frac{1}{1-z^{w_i}}.
\end{equation}
In particular, acyclic orientations with different multisets of path
weights give nonisomorphic graded endomorphism algebras.  For type $D_4$,
orienting the three edges toward the central vertex gives weights
$1,1,1,4$, whereas orienting them away from the central vertex gives
$1,2,2,2$.
\end{corollary}

\begin{proof}
The filtration in \eqref{eq:directed-path-classification} is the weighted
degree filtration on a polynomial algebra.  Its associated graded algebra
is therefore polynomial on the initial forms of the variables, with the
displayed degrees.  Formula \eqref{eq:path-hilbert-series} follows.  A
graded algebra isomorphism preserves the Hilbert series.  This series also
determines the multiset of positive weights.  If $H(z)$ denotes the right
side of \eqref{eq:path-hilbert-series}, then in the formal power-series
ring $k[[z]]$, which has characteristic zero under the standing
hypotheses, its logarithmic derivative is
\[
 z\frac{d}{dz}\log H(z)
 =\sum_{n\geq1}\left(\sum_{w_i\mid n}w_i\right)z^n.
\]
Write $c_n=\sum_{w_i\mid n}w_i$ and let $m_n$ be the number of indices
with $w_i=n$.  The recursion
\[
 n m_n=c_n-\sum_{d\mid n,\ d<n}d m_d
\]
recovers every $m_n$, starting with $n=1$.  Hence different multisets of
weights have different Hilbert series.  The stated $D_4$ weights follow
directly from \eqref{eq:path-weights}.
\end{proof}

\begin{corollary}
\label{cor:A2-quantum-plane}
In type $A_2$, with the order $1<2$, the module above is
\[
 K_1v_{m,n}=q^{2m-n}v_{m,n},\qquad
 K_2v_{m,n}=q^{-m+2n}v_{m,n},
\]
\[
 E_1v_{m,n}=v_{m+1,n},\qquad
 E_2v_{m,n}=q^m v_{m,n+1}.
\]
Its filtered endomorphism algebra is $k[X,Y]$ with
$\deg X=1$ and $\deg Y=2$.
\end{corollary}

The path weights depend on the chosen acyclic orientation.  For the
linearly oriented diagram of type $A_\ell$ they are
$1,2,\ldots,\ell$.  More generally, formula \eqref{eq:path-weights}
measures the propagation of a defect against the orientation: applying a
later root defect increases the number of earlier defects that can remain
nonzero.  This gives a combinatorial invariant of the chosen
quantum-space quotient rather than an invariant of the unoriented Dynkin
diagram alone.

\section{Bar cochains and filtered composition}

We recall the bar-resolution realization of the defect complex and then
construct its cochain composition.  The Ext comparison itself is standard.
The explicit Hopf Alexander--Whitney formula below realizes Yoneda
composition on cohomology.  The filtration in this section is exclusively
the absolute augmentation filtration $\cD_H^r$ from Section~3.  It is not
the relative filtration $\cD_{H\mid K}^r$ of Section~5, and the results of
this section do not assert a cohomological realization of the relative
principal symbols.
Let $V$ be a left $H$-module.  Give it the right $H$-action
$v\cdot x=\varepsilon(x)v$.  The Hochschild cochain complex of $H$ with
coefficients in this bimodule is
\[
 C^n(H,V)=\kHom(H^{\otimes n},V),\qquad n\ge0,
\]
with differential
\begin{align}
 (dc)(x_1,\ldots,x_{n+1})={}&x_1c(x_2,\ldots,x_{n+1})\notag\\
 &+\sum_{i=1}^{n}(-1)^ic(x_1,\ldots,x_ix_{i+1},\ldots,x_{n+1})\notag\\
 &+(-1)^{n+1}c(x_1,\ldots,x_n)\varepsilon(x_{n+1}).
 \label{eq:differential}
\end{align}
In degree zero, $(dv)(x)=xv-\varepsilon(x)v=\langle x,v\rangle$.

For modules $M,N$, abbreviate
\[
 C^\bullet(H;M,N)=C^\bullet(H,\kHom(M,N)).
\]

\begin{theorem}\label{thm:ext}
For arbitrary left $H$-modules $M,N$, there is a natural isomorphism
\begin{equation}\label{eq:ext}
 H^n(C^\bullet(H;M,N))\cong\Ext_H^n(M,N),\qquad n\ge0.
\end{equation}
In degree zero this identifies cocycles with ordinary intertwiners.
No finite-dimensionality assumption is required.
\end{theorem}

\begin{proof}
Let $\Barres_\bullet(H)$ be the unnormalized left bar resolution of the
trivial module $k$:
\[
 \Barres_n(H)=H^{\otimes(n+1)}.
\]
Each term is free as a left $H$-module through multiplication on the first
factor.  For $n\geq1$, its boundary is
\begin{align*}
 d(h_0\otimes\cdots\otimes h_n)
 ={}&\sum_{i=0}^{n-1}(-1)^i
 h_0\otimes\cdots\otimes h_ih_{i+1}\otimes\cdots\otimes h_n\\
 &+(-1)^n h_0\otimes\cdots\otimes h_{n-1}\varepsilon(h_n),
\end{align*}
and the augmentation is $\varepsilon:\Barres_0(H)=H\to k$.  The standard
contracting homotopy of the augmented complex inserts $1$ in the first
tensor position.  Therefore $\Barres_\bullet(H)\otimes M$, with diagonal
$H$-action, is a resolution of $M$.  Its augmentation is the $H$-linear map
\[
 \varepsilon_M:H\otimes M\longrightarrow M,
 \qquad h\otimes m\longmapsto\varepsilon(h)m,
\]
and the bar contracting homotopy tensored with $\id_M$ proves exactness.
This contracting homotopy is only $k$-linear; it is used to prove
exactness, not to split the resolution in the category of $H$-modules.

We verify projectivity explicitly.  Define
\begin{align*}
 \Phi_n:H\otimes H^{\otimes n}\otimes M&\longrightarrow
 H^{\otimes(n+1)}\otimes M,\\
 \Phi_n(h\otimes\mathbf x\otimes m)
 &=\sum h_{(1)}\otimes x_1\otimes\cdots\otimes x_n
       \otimes h_{(2)}m,
\end{align*}
where the source has the regular action on the first factor and trivial
action on the remaining vector space.  Its inverse is
\[
 h_0\otimes\mathbf x\otimes m\longmapsto
 \sum h_{0(1)}\otimes\mathbf x\otimes S(h_{0(2)})m.
\]
The antipode identities show that these maps are inverse and $H$-linear.
Thus every term is free, even when $M$ is infinite dimensional.

By Proposition~\ref{prop:adjunction}, degreewise currying gives
\begin{align*}
 \operatorname{Hom}_H(\Barres_n(H)\otimes M,N)
 &\cong\operatorname{Hom}_H(\Barres_n(H),\kHom(M,N))\\
 &\cong\kHom(H^{\otimes n},\kHom(M,N)).
\end{align*}
A direct substitution of the bar boundary shows that the induced
coboundary is exactly \eqref{eq:differential}.  The cohomology of the left
side is $\Ext_H^n(M,N)$, proving \eqref{eq:ext}.  For $n=0$, the
differential is the equivariance defect, and its kernel is
$\operatorname{Hom}_H(M,N)$.
\end{proof}

For later use, recall that a cohomology class in
$\operatorname{Hom}_H(P(M),N)$ may be lifted along a projective resolution
$P(N)\to N$ to a chain map $P(M)\to P(N)[n]$.  Composition of such lifts
represents the Yoneda product.  This standard description will be used
below only to identify the product induced by the explicit bar diagonal
\cite[Chapter~VIII]{CartanEilenberg}.

\begin{lemma}\label{lem:hopf-aw}
Let $\Barres_\bullet(H)$ be the bar resolution in
Theorem~\ref{thm:ext}.  For $p+q=n$, let
\[
 D_{p,q}:\Barres_n(H)\longrightarrow
 \Barres_p(H)\otimes\Barres_q(H)
\]
be the $H$-linear map determined by
\begin{align}
 1\otimes x_1\otimes\cdots\otimes x_n
 \longmapsto
 \sum{}&
 (1\otimes x_{1(1)}\otimes\cdots\otimes x_{p(1)})\notag\\
 &\otimes
 (x_{1(2)}\cdots x_{p(2)}
   \otimes x_{p+1}\otimes\cdots\otimes x_n).
 \label{eq:AW-Hopf}
\end{align}
Then $D_n=\sum_{p+q=n}D_{p,q}$ defines a counital, strictly
coassociative chain map
\[
 D:\Barres_\bullet(H)\longrightarrow
 \Barres_\bullet(H)\otimes\Barres_\bullet(H),
\]
where the tensor product has its usual total differential and the diagonal
$H$-action.
\end{lemma}

\begin{proof}
The bar resolution is the chain complex associated with the simplicial bar
construction.  Applying the Hopf coproduct to each bar entry gives a
simplicial map from this construction to its degreewise tensor square.
Composing it with the ordinary Alexander--Whitney map gives $D$.  On the
generator $1\otimes x_1\otimes\cdots\otimes x_n$, its $(p,q)$ component
is exactly \eqref{eq:AW-Hopf}; the product
$x_{1(2)}\cdots x_{p(2)}$ appears when the leading terms in the second
bar factor are multiplied.

For completeness, the chain-map identity follows from the simplicial
identities.  Away from the cut between $p$ and $p+1$, the relevant face
multiplies two adjacent entries, and the equality on the two sides is
$\Delta(x_ix_{i+1})=\Delta(x_i)\Delta(x_{i+1})$.  At the two outer faces
it is the counit identity.  The two terms at the cut are the two summands
of the total differential and have the required sign.  Thus
\[
 d_{\mathrm{tot}}D=Dd.
\]
Naturality and strict coassociativity of the ordinary Alexander--Whitney
map, together with $(\Delta\otimes\id)\Delta=(\id\otimes\Delta)\Delta$,
give
\[
 (D\otimes\id)D=(\id\otimes D)D.
\]
The two counit identities for $D$ follow from the augmentation of the bar
resolution.  Finally, the formula is $H$-linear because the first tensor
factor in \eqref{eq:AW-Hopf} carries the regular action and the target has
the diagonal action.
\end{proof}

\begin{theorem}\label{thm:cup-composition}
For
\[
 b\in C^p(H;N,L),\qquad c\in C^q(H;M,N),
\]
define $b\smile c\in C^{p+q}(H;M,L)$ by
\begin{align}
 &(b\smile c)(x_1,\ldots,x_{p+q}) \notag\\
 &\quad =
 \sum b(x_{1(1)},\ldots,x_{p(1)})\circ
 \Bigl((x_{1(2)}\cdots x_{p(2)})\cdot
 c(x_{p+1},\ldots,x_{p+q})\Bigr).
 \label{eq:cup-composition}
\end{align}
For $p=0$, the empty product $x_{1(2)}\cdots x_{p(2)}$ is understood to
be $1$.  Then
\begin{equation}\label{eq:cup-leibniz}
 d(b\smile c)=(db)\smile c+(-1)^p b\smile(dc).
\end{equation}
The product is associative, its degree-zero units are the identity
intertwiners, and the product induced on cohomology agrees under
\eqref{eq:ext} with Yoneda composition.
\end{theorem}

\begin{proof}
Apply the cochains represented by $b$ and $c$ to the two tensor factors
of the diagonal in Lemma~\ref{lem:hopf-aw}, and then apply the composition
pairing
\[
 \kHom(N,L)\otimes\kHom(M,N)\longrightarrow\kHom(M,L).
\]
This pairing is $H$-linear by Proposition~\ref{prop:composition}.
The resulting expression is exactly \eqref{eq:cup-composition}.
Since the Alexander--Whitney map is a chain map, applying the cochain
differential to the first tensor factor gives $(db)\smile c$, while
applying it to the second factor crosses a cochain of degree $p$ and
therefore gives $(-1)^p b\smile(dc)$.  This proves
\eqref{eq:cup-leibniz}.
Coassociativity of the Alexander--Whitney diagonal and associativity of
composition give associativity of $\smile$.  The augmentation supplies
the degree-zero units.

To identify the product on cohomology, let $c$ be a cocycle.  The component
of $D$ with second degree $q$, followed by $c$ on that factor, gives the
standard comparison lift from the bar resolution of $M$ to the shifted bar
resolution of $N$.  Applying the analogous lift for $b$ gives the cochain
$b\smile c$.  Under the usual correspondence between comparison maps and
extensions, composition of these lifts represents the splice of the
extensions represented by $c$ and $b$.  This is Yoneda composition
\cite[Chapter~VIII]{CartanEilenberg}.
\end{proof}

Formula \eqref{eq:cup-composition} explains why simply dividing the
arguments of $b$ and $c$ is insufficient.  The first $p$ arguments must
also be split by the coproduct, and their second tensor factors act on
the value of $c$.  For a group algebra, the formula reduces to the
standard inhomogeneous group-cohomology cup product.

\subsection{Normalized cochains and the absolute augmentation filtration}

On the unnormalized complex put
\[
 F_r C^n(H;M,N)=
 \kHom(H^{\otimes n},\cD^r_H(M,N)).
\]
The normalized complex is
\[
 \overline C^n(H;M,N)=\kHom(I^{\otimes n},\kHom(M,N)).
\]
It is quasi-isomorphic to the unnormalized complex.  The filtration on
coefficients induces
\[
 F_r\overline C^n(H;M,N)=
 \kHom(I^{\otimes n},\cD^r_H(M,N)).
\]

\begin{theorem}\label{thm:filtered-complex}
The differential preserves the increasing filtration $F_r$ on both the
unnormalized and normalized complexes.  For all $a,b,p,q\geq0$,
\begin{align}
 F_a C^p(H;N,L)\smile F_b C^q(H;M,N)
 &\subseteq F_{a+b}C^{p+q}(H;M,L), \label{eq:filtered-cup}\\
 F_a\overline C^p(H;N,L)\smile F_b\overline C^q(H;M,N)
 &\subseteq F_{a+b}\overline C^{p+q}(H;M,L). \notag
\end{align}
Consequently, the unions
\[
 F_\infty\overline C^\bullet(H;M,N)
 =\bigcup_{r\geq0}F_r\overline C^\bullet(H;M,N)
\]
form a filtered differential graded category.  For $R\geq0$, let
$E_t^{i,j}(R;M,N)$ be the cohomological spectral sequence of the finite
increasing filtration
\[
 0=F_{-1}\overline C^\bullet\subset F_0\overline C^\bullet
 \subset\cdots\subset F_R\overline C^\bullet.
\]
It has
\[
 E_0^{i,j}(R;M,N)=
 \frac{F_i\overline C^{,i+j}(H;M,N)}
      {F_{i-1}\overline C^{,i+j}(H;M,N)},
 \qquad 0\leq i\leq R,
\]
and converges to $H^{i+j}(F_R\overline C^\bullet(H;M,N))$.  For finite
$R,S$ and every page $t\geq0$, cup-composition induces natural pairings
\[
 E_t^{i,j}(R;N,L)\otimes E_t^{i',j'}(S;M,N)
 \longrightarrow
 E_t^{i+i',j+j'}(R+S;M,L).
\]
The pairing on the abutments is induced by $\smile$ and, after inclusion
in the full cochain complexes, agrees with Yoneda composition.  On the
$E_0$ page, coefficient composition is the absolute graded composition of
Corollary~\ref{cor:graded}.  The associated graded category for the relative
filtration of Section~5 is not involved in this assertion.  The order-zero coefficient piece is
$\operatorname{Hom}_H(M,N)$.
\end{theorem}

\begin{proof}
The subspace $\cD^r_H(M,N)$ is an $H$-submodule.  Indeed, if $a\in
I^{r+1}$, then $a\cdot(x\cdot f)=(ax)\cdot f=0$ because $I^{r+1}$ is a
right ideal.  Every term of \eqref{eq:differential} therefore remains in
$\cD^r_H(M,N)$, which proves preservation by the differential.

Let $b$ take values in $\cD^a_H(N,L)$ and let $c$ take values in
$\cD^b_H(M,N)$.  Since $\cD^b_H(M,N)$ is an $H$-submodule, every value
of
\[
 (x_{1(2)}\cdots x_{p(2)})\cdot
 c(x_{p+1},\ldots,x_{p+q})
\]
still lies in $\cD^b_H(M,N)$.  Theorem~\ref{thm:main} now places each
summand of \eqref{eq:cup-composition} in $\cD^{a+b}_H(M,L)$ and proves
the first inclusion in \eqref{eq:filtered-cup}.  The
Alexander--Whitney diagonal descends to the normalized bar resolution,
so the same proof gives the normalized inclusion.

The differential graded category assertion follows from
\eqref{eq:cup-leibniz}, associativity, and \eqref{eq:filtered-cup}.
For fixed $R$, the filtration of $F_R\overline C^\bullet$ is finite, so
its spectral sequence converges to the cohomology of that subcomplex.  Its
zeroth page is the displayed associated graded complex.  The filtered
pairing \eqref{eq:filtered-cup} has filtration degree zero after the
source filtrations are given the product filtration.  The usual
cycle-and-boundary construction therefore induces the stated pairing on
every page, and the Leibniz identity makes each differential a derivation
for this pairing.  The abutment pairing is induced by $\smile$, hence maps
to Yoneda composition by Theorem~\ref{thm:cup-composition}.  Passing
\eqref{eq:filtered-cup} to the $E_0$ page gives the graded absolute
coefficient composition of Corollary~\ref{cor:graded}.  Finally,
$\cD^0_H(M,N)=\operatorname{Hom}_H(M,N)$ by
\eqref{eq:order-zero}.
\end{proof}

\subsection{Adjoint coefficients and the derived center}

Let $H_{\ad}$ denote $H$ with the left adjoint action
\[
 x\triangleright y=\sum_{(x)}x_{(1)}yS(x_{(2)}).
\]
The associated defect is the reduced adjoint bracket
\begin{equation}\label{eq:adj-bracket}
 \langle x,y\rangle=\sum_{(x)}x_{(1)}yS(x_{(2)})-\varepsilon(x)y.
\end{equation}

\begin{proposition}\label{prop:center}
The order-zero adjoint invariants are the center:
\[
 \cD^0_H(k,H_{\ad})=(H_{\ad})^H=Z(H).
\]
More generally,
\[
 \cD^r_H(k,H_{\ad})=\{y\in H:I^{r+1}\triangleright y=0\}
\]
is a multiplicative filtration of the subalgebra
$\cD^\infty_H(k,H_{\ad})\subseteq H$.
\end{proposition}

\begin{proof}
An element $y$ is adjoint invariant precisely when
$\sum x_{(1)}yS(x_{(2)})=\varepsilon(x)y$ for every $x$.  In that case,
\[
 xy=\sum (x_{(1)}\triangleright y)x_{(2)}
    =\sum\varepsilon(x_{(1)})yx_{(2)}=yx.
\]
Conversely, if $y$ is central, then
\[
 x\triangleright y
 =\sum x_{(1)}yS(x_{(2)})
 =y\sum x_{(1)}S(x_{(2)})
 =\varepsilon(x)y.
\]
The second statement is the definition.  Since the adjoint action makes $H$ an
$H$-module algebra, one has
$x\triangleright(yz)=\sum(x_{(1)}\triangleright y)
(x_{(2)}\triangleright z)$.  Lemma~\ref{lem:filtered-coproduct} and the
same pigeonhole argument as in Theorem~\ref{thm:main} prove
multiplicativity.
\end{proof}

We now recall the standard comparison between Hochschild cohomology and
Hopf cohomology with adjoint coefficients.  On a Hochschild cochain
$f:H^{\otimes n}\to H$ we use
\begin{align*}
 (d_{\mathrm{Hoch}}f)(x_1,\ldots,x_{n+1})={}&
 x_1f(x_2,\ldots,x_{n+1})\\
 &+\sum_{i=1}^{n}(-1)^if(x_1,\ldots,x_ix_{i+1},\ldots,x_{n+1})\\
 &+(-1)^{n+1}f(x_1,\ldots,x_n)x_{n+1}.
\end{align*}

\begin{theorem}\label{thm:center}
There is a natural isomorphism of graded vector spaces
\begin{equation}\label{eq:center-comparison}
 \HH^\bullet(H,H)\cong\Ext_H^\bullet(k,H_{\ad}).
\end{equation}
When the usual cup products are used, this is an isomorphism of graded
algebras.  In degree zero it is the identification
$Z(H)=(H_{\ad})^H$ of Proposition~\ref{prop:center}.
\end{theorem}

\begin{proof}
Regard $H$ as an $H$-bimodule and identify $H^e=H\otimes H^{\mathrm{op}}$.
Restriction along
\[
 H\longrightarrow H^e,\qquad
 x\longmapsto\sum x_{(1)}\otimes S(x_{(2)})
\]
gives the adjoint action.  On cochains, define
\begin{align*}
 \Phi(f)(x_1,\ldots,x_n)
 &=\sum f(x_{1(1)},\ldots,x_{n(1)})
       S(x_{1(2)}\cdots x_{n(2)}),\\
 \Psi(g)(x_1,\ldots,x_n)
 &=\sum g(x_{1(1)},\ldots,x_{n(1)})
       x_{1(2)}\cdots x_{n(2)}.
\end{align*}
Here $f$ is a Hochschild cochain with values in the regular bimodule $H$,
and $g$ is a Hopf cochain with values in $H_{\ad}$.  The antipode
identities show that $\Phi\Psi$ and $\Psi\Phi$ are identity maps.  For
example, the rightmost Hochschild face, which multiplies the value of a
cochain by $x_{n+1}$, becomes after applying $\Phi$ the counit term
$g(x_1,\ldots,x_n)\varepsilon(x_{n+1})$.  The leftmost face becomes the
adjoint action, and the internal faces are unchanged because $\Delta$ is
an algebra map.  Thus $\Phi$ and $\Psi$ commute with the differentials.

The coproduct factors in these formulas also identify the
Alexander--Whitney diagonals on the two bar constructions.  Consequently
the cup products correspond.  See \cite[Section 3]{LachowskaQi} for the
multiplicative comparison in the finite-dimensional setting; the formulas
above do not require finite-dimensionality.
\end{proof}

For a small quantum group, the full Hochschild cohomology
is its derived center.  Geometric descriptions identify important blocks
of this algebra with cohomology of polyvector fields on the Springer
resolution \cite{BezrukavnikovLachowska}, while explicit computations of
centers in low rank are carried out in \cite{LachowskaQiCenter}.  The
absolute filtration in Proposition~\ref{prop:center} supplies an elementary
augmentation-theoretic filtration on the degree-zero part and, through
Theorem~\ref{thm:filtered-complex}, on a cochain model for higher degrees.
This statement is separate from the relative quantum-Borel filtration of
Section~5.

\section{Lie-de Rham cocycles and quantum screening applications}

Ginzburg and Schechtman construct total Lie-de Rham cocycles and the
associated Ext-valued pairings
\cite[Theorem~3.12 and Corollary~3.13]{GinzburgSchechtman}.
We recall the algebraic pairing in the notation of Section~6 to distinguish
its cohomological degree from the order of a higher intertwiner.

Let $\Omega^\bullet$ be a complex of vector spaces concentrated in degrees
$0,\ldots,a$, with trivial $H$-action.  Set
\[
 C^{p,q}=C^p(H;M,N\otimes\Omega^q),
 \qquad D=d_H+(-1)^p d_\Omega\quad\text{on }C^{p,q},
\]
where $d_H$ is \eqref{eq:differential} and $d_\Omega$ is induced by the
coefficient differential.  These two differentials commute before the
indicated totalization sign.  Suppose
\[
 \mathcal V=\sum_{p=0}^a\mathcal V^{p,a-p},
 \qquad \mathcal V^{p,a-p}\in C^{p,a-p},
 \qquad D\mathcal V=0.
\]
For a degree-$m$ functional $\gamma:\Omega^m\to k$ satisfying
$\gamma d_\Omega=0$, evaluation defines
\[
 c_\gamma=(\id_N\otimes\gamma)\mathcal V^{a-m,m}
 \in C^{a-m}(H;M,N).
\]
The component equations of $D\mathcal V=0$ give $d_Hc_\gamma=0$.
If $\gamma$ is changed by a boundary in the dual complex, the preceding
component of $\mathcal V$ shows that $c_\gamma$ changes by a
$d_H$-coboundary.  Thus Theorem~\ref{thm:ext} gives
\[
 H_m\bigl((\Omega^\bullet)^\vee\bigr)
 \longrightarrow\Ext_H^{a-m}(M,N),
 \qquad [\gamma]\longmapsto[c_\gamma].
\]
Here $(\Omega^\bullet)^\vee$ is the chain complex with degree-$m$ term
$\operatorname{Hom}_k(\Omega^m,k)$ and boundary given by precomposition
with $d_\Omega$.  This is the standard total-cocycle pairing; it makes no
claim that the resulting Ext classes are nonzero.  The tensor products
in this explanation are algebraic.  Formal-series applications require
the corresponding completed coefficient spaces and compatible evaluation
pairings, not an identification of
$\operatorname{Hom}_k(M,N\otimes\Omega^q)$ with
$\kHom(M,N)\otimes\Omega^q$.

For total degree one, write $\mathcal V=\alpha+\beta$, with
$\alpha\in C^{0,1}$ and $\beta\in C^{1,0}$.  The equations are
\[
 d_\Omega\alpha=0,\qquad
 d_H\alpha=d_\Omega\beta,\qquad d_H\beta=0.
\]
The middle equation expresses equivariance up to a coefficient
differential.  By itself it does not imply the last equation.  In
particular, the fact that a defect $d_Hf$ is a cochain does not make it a
nonzero Ext class: in the full coefficient complex it is a coboundary.

In \cite{Sebbar1998}, operator-valued $q$-differential forms between Verma
modules satisfy the analogous two families of equations.  The
$q$-de Rham differential controls the configuration variables, while
\eqref{eq:differential} controls equivariance.  For
$U_q(\widehat{\mathfrak{sl}}_2)$, quantum screening operators provide the
canonical cocycles \cite{SebbarCMP}.  Pairing with admissible twisted
$q$-de Rham cycles gives maps of the form
\[
 H_m^{q\text{-}\mathrm{dR}}(\mathcal L)
 \longrightarrow
 \Ext_{U_q(\mathfrak g)}^{a-m}(M(\lambda'),M(\lambda)).
\]
The bar model in Theorem~\ref{thm:ext} supplies a common algebraic
framework for the target, including for infinite-dimensional modules.

Cohomological degree $a-m$ is independent of the order defined by
augmentation-ideal powers.  Proposition~\ref{prop:indecomp-Ug} and
Theorem~\ref{thm:quantum-collapse} give cases in which all finite absolute
orders reduce to ordinary intertwiners; this does not imply the vanishing
of higher Ext groups.  The relative construction in Section~5 instead
requires Cartan equivariance and measures root defects by powers of the
projection kernel.  Its symbol obstruction arises from the specific
extension $0\to C_r\to A_r\to A_{r-1}\to0$, not from a de Rham
integration pairing.  No identification of these two constructions, and
no relative cochain filtration on the screening complex, is asserted.
Likewise, free-field bosonization in the screening literature should not
be confused with the Radford--Majid biproduct $R\#K$ used in
Theorem~\ref{thm:bosonization}.

The wider geometric setting begins with the relation between hyperplane
arrangements and Lie algebra homology developed by Schechtman and
Varchenko \cite{SchechtmanVarchenko1991}.  Their later work
\cite{SchechtmanVarchenko2017} constructs a map from a twisted de Rham
complex on a punctured projective line to a Lie algebra chain complex
with affine Verma-module coefficients.  The derived KZ comparison
\cite[Theorem~3.8 and Corollary~3.9]{SchechtmanVarchenko2022} identifies
the induced KZ and Gauss--Manin connections in higher degrees, with
isomorphisms for generic parameter.  Its filtration is in the base
configuration variables; it is not the augmentation-order filtration
used here.

Two related geometric developments are the realization of the
$\mathfrak{sl}_2$ BGG resolution by Aomoto complexes in
\cite{FalkSchechtmanVarchenko}, and the configuration-space construction
of representations of small quantum groups at roots of unity in
\cite{BezrukavnikovFinkelbergSchechtman}.  The latter concerns a different
parameter regime from the generic-$q$ Borel calculations of Section~5.
These works provide geometric context without supplying, or being
identified with, the relative symbol calculus of this paper.

\subsection*{Acknowledgements.} This research was partially supported by an NSERC Discovery Grant.

\end{document}